\documentclass[10pt]{amsart}
\usepackage[T1]{fontenc}
\usepackage{amsmath,amssymb,mathtools,mathrsfs,microtype}
\usepackage[margin=0.95in]{geometry}
\usepackage{tikz-cd,enumitem}
\usepackage[colorlinks=true,citecolor=blue,linkcolor=blue,urlcolor=blue]{hyperref}
\allowdisplaybreaks[2]
\newtheorem{theorem}{Theorem}[section]
\newtheorem{proposition}[theorem]{Proposition}
\newtheorem{lemma}[theorem]{Lemma}
\newtheorem{corollary}[theorem]{Corollary}
\theoremstyle{definition}
\newtheorem{definition}[theorem]{Definition}
\theoremstyle{remark}
\newtheorem{remark}[theorem]{Remark}
\theoremstyle{plain}
\newtheorem*{theoremA}{Theorem A}
\newcommand{\Q}{\mathbb Q}
\newcommand{\Z}{\mathbb Z}
\newcommand{\C}{\mathbb C}
\newcommand{\Gm}{\mathbb G_m}
\newcommand{\A}{\mathbb A}
\newcommand{\cA}{\mathcal A}
\newcommand{\cU}{\mathcal U}
\newcommand{\cM}{\mathfrak M}
\newcommand{\HH}{\operatorname{HH}}
\newcommand{\DM}{\mathrm{DM}_{\mathrm{gm}}}
\newcommand{\Spec}{\operatorname{Spec}}
\newcommand{\Hom}{\operatorname{Hom}}

\newcommand{\id}{\operatorname{id}}

\newcommand{\Gr}{\operatorname{Gr}}
\newcommand{\age}{\operatorname{age}}
\newcommand{\sat}{\operatorname{Sat}}
\newcommand{\CH}{\operatorname{CH}}
\newcommand{\et}{\'{e}tale}
\newcommand{\Morb}{M_{\mathrm{orb}}}
\title[Motivic mirror symmetry for sink-recurrent cluster varieties]
{Motivic mirror symmetry\\for sink-recurrent cluster varieties}
\author{Haoyang Liu}
\address{Department of Mathematics, University of California, Santa Barbara,
CA 93106, USA}
\email{haoyangliu@ucsb.edu}
\author{Tianle Liu}
\address{University of Southern California}
\email{tianleli@usc.edu}
\date{}
\hypersetup{pdftitle={Motivic mirror symmetry for sink-recurrent cluster varieties},pdfauthor={}}
\begin{document}
\begin{abstract}
We prove motivic mirror symmetry for full-rank sink-recurrent cluster
varieties. We allow skew-symmetrizable principal parts and invert all
frozen variables. The rational geometric motive of such a variety is
isomorphic to the additive orbifold motive of the mirror stack
of Hu and Shende. This gives comparisons of cohomology,
higher Chow groups, and algebraic $K$-theory with rational coefficients.
The proof combines quantum deformation with the splitting of these
motives into Tate summands. For really full-rank skew-symmetric seeds,
a finite \et\ morphism realizes the comparison. We also study this
geometric comparison in families over the projective line.
\end{abstract}
\maketitle

\section{Introduction}

Fock and Goncharov introduced dual $\mathcal A$ and $\mathcal X$
cluster varieties and formulated duality conjectures relating canonical
bases of regular functions to integral tropical points of the
dual variety~\cite{FG}. Arg\"uz and Bousseau identified
Gross--Siebert mirrors of cluster varieties with degenerations of their
Fock--Goncharov duals~\cite{ArguzBousseau}.

On the symplectic side, Gammage and Le proved homological mirror symmetry
for truncated cluster varieties~\cite{GammageLe}. Hu and Shende prove
homological mirror symmetry for full-rank sink-recurrent cluster
varieties~\cite{HS}.

In this paper, we prove a motivic version of mirror symmetry for
full-rank sink-recurrent seeds with skew-symmetrizable principal part.
Our main theorem identifies the motive of a cluster variety
with the orbifold motive of its mirror. It yields comparisons of their
cohomology, higher Chow groups, and algebraic $K$-theory with rational
coefficients.

Let
\[
 A=\begin{pmatrix}B\\ C\end{pmatrix}\in M_{N\times n}(\Z),
 \qquad N=n+f,
\]
have full column rank and sink-recurrent principal part, and let
$D=\operatorname{diag}(d_1,\ldots,d_n)$
be a positive integral matrix such that $DB$ is skew-symmetric.
We take $D=I_n$ when $B$ is skew-symmetric.
Let $\cA(A)$ be the ordinary cluster algebra
of geometric type, with all frozen variables inverted, and write
$X(A)=\Spec\cA(A)$ for the associated cluster variety.
Sink recurrence means that a finite sequence of mutations and recursive
sink freezings reduces the seed to isolated seeds
(Definition~\ref{geom:sr}). It implies that $X(A)$ is smooth and of
finite type.

For an $n\times n$ exchange matrix $B_0$, let $Y_{B_0}$ be the cluster variety defined by the
extended exchange matrix $\bigl(\begin{smallmatrix}B_0\\I_n\end{smallmatrix}\bigr)$;
the bottom block gives its $n$ principal coefficients. Following
Hu and Shende~\cite{HS}, define the mirror to be the quotient stack
\begin{equation}\label{intro:mirror}
 P_A=Y_{-B^{\mathsf T}}\times\Gm^f,
 \qquad \cM_A=[P_A/H],\qquad H=\Gm^n.
\end{equation}
Write $X_1,\ldots,X_n$ for the initial mutable cluster variables on
$Y_{-B^{\mathsf T}}$, and $Z_1,\ldots,Z_n$ for its frozen variables.
Let $U_1,\ldots,U_f$ be the coordinates on the factor $\Gm^f$.
For $q=(q_1,\ldots,q_n)\in H$, the action is
\begin{equation}\label{intro:action}
 q\cdot X_i=q_iX_i,\qquad
 q\cdot Z_i=\Bigl(\prod_{j=1}^n q_j^{b_{ij}}\Bigr)Z_i,\qquad
 q\cdot U_\ell=\Bigl(\prod_{j=1}^n q_j^{c_{\ell j}}\Bigr)U_\ell.
\end{equation}
For sink-recurrent $A$, this quotient contains the Langlands-dual
Fock--Goncharov $\mathcal X$-variety as an open substack
(Proposition~\ref{geom:fg-affinization}).
The Langlands-dual matrix entering its finite quotient presentation is
\begin{equation}\label{intro:dual}
 A^\vee=\begin{pmatrix}-B^{\mathsf T}\\C\operatorname{diag}(p_i)\end{pmatrix},
 \qquad p_i=\frac{\operatorname{lcm}(d_1,\ldots,d_n)}{d_i}.
\end{equation}
The least common multiple of an empty set is taken to be $1$.
The matrix $A^\vee$ is again full-rank and sink-recurrent.
Write $X^\vee=X(A^\vee)$.
Proposition~\ref{geom:mirror} presents $\cM_A$ as $[X^\vee/G_0]$
for an explicit finite diagonal group $G_0$.
When $n=0$, both $X(A)$ and $\cM_A$ are $\Gm^f$ and the
motivic comparison is immediate.

We use rational geometric motives with
$M(\Gm)=\Q\oplus\Q(1)[1]$, and write $M\{a\}=M(a)[2a]$.
For a smooth quotient by a finite abelian group $G$ whose fixed
sectors have integral ages,
its additive orbifold motive is
\begin{equation}\label{intro:orbifold}
 \Morb([X/G])=
 \bigoplus_{g\in G}\ \bigoplus_a M(X^g_a)^G\{a\}.
\end{equation}
Here $X^g_a$ is the union of the connected components of age $a$;
group-permuted components are combined before taking invariants.
The invariant motive is the rational motive of the corresponding
finite quotient stack~\cite[Lemma~3.1]{Choudhury}.
The finite quotient presentation in Proposition~\ref{geom:mirror}
makes~\eqref{intro:orbifold} applicable to $\cM_A$. Its ages are
half the codimensions of the fixed components, by the symplectic
GSV form~\cite{GSV2005} in even dimension and by adjoining one inert torus factor
in odd dimension.

\begin{theoremA}
Let $A$ be a full-rank sink-recurrent extended exchange matrix with
skew-symmetrizable principal part. Then there is an isomorphism
\[
 M(X(A))\simeq\Morb(\cM_A)
 \qquad\text{in }\DM(\C;\Q).
\]
\end{theoremA}

Theorem~\ref{main:character} refines this comparison character by
character, pairing each component character with the fixed sector
specified by the integral exchange lattice and its age shift. The
additive isomorphism depends on choices of bases in the Tate
multiplicity spaces. We also obtain the
compactly supported comparison and, after products with smooth
varieties, isomorphisms of higher Chow groups and rational algebraic
$K$-theory (Corollary~\ref{main:chow-k}).

When $A^{\mathsf T}:\Z^N\to\Z^n$ is surjective, we call $A$
\emph{really full-rank}. In this case the mirror is the spectrum of
another cluster algebra. When $B$ is skew-symmetric, an explicit finite
\et\ morphism $p_A:X(A)\to\cM_A$ realizes the comparison.
Its extension to deformation spaces of torus-stable smooth
closed pairs gives a geometric interpretation over $\mathbb P^1$
(Section~\ref{sec:p1}).

\subsection{Proof strategy}

Write $X=X(A)$. Sink covers reduce the source, the dual variety, and
their diagonal fixed loci to isolated charts. Their stratifications
give mixed Tate motives, and curious Lefschetz gives split mixed Tate
Hodge structures. An arithmetic argument over a cyclotomic number field
then splits the motives into shifted Tate summands
(Theorem~\ref{geom:splitting}). We compare their multiplicities using
two formal realizations of one quantum cluster algebra and its
Kodaira--Spencer action.

Suppose first that $N=2d$. An integral nondegenerate alternating
completion of $A\operatorname{diag}(p_i)$ defines a compatible quantum cluster
algebra $\cA_t$, with $\cA_1=\mathcal O(X)$; write
$\cA_{\mathrm{gen}}=\cA_t\otimes_{\C[t^{\pm1}]}\C(t)$.
For a root of unity $\rho$ of sufficiently divisible order and a suitable
finite diagonal group $Q$, the crossed product
$\mathcal B_t=\cA_t\rtimes Q$ has
$Z(\mathcal B_\rho)=\mathcal O(X^\vee)$. The central coordinates are
powers of quantum cluster variables multiplied by group elements.
Brown--Gordon's Poisson-order theorem~\cite[Theorem~4.2]{BrownGordon}
then makes $\mathcal B_\rho$ Azumaya. For $G=Q^\vee$,
crossed-product duality gives
\[
 \mathcal B_t\rtimes G\cong
 \operatorname{Mat}_{|Q|}(\cA_t).
\]
It pairs source characters with fixed sectors of $[X^\vee/G]$,
whose invariant motives agree with those of the mirror presentation
$[X^\vee/G_0]$ (Proposition~\ref{geom:inertia-comparison}).

Theorem~\ref{ql:completion} compares generic Hochschild homology with
its continuous versions at $t=e^\hbar$ and $t=\rho e^\hbar$.
At isolated seeds, contractions with uniformly bounded formal poles
control the growing Laurent supports allowed by completion. Ore
localizations and sink covers extend the comparison to the whole
sink-recurrent variety.

For a source character $\chi$ and its paired element $g\in G$,
symplectic deformation quantization gives two realizations of
$\HH_i(\cA_{\mathrm{gen}})_\chi$:
\[
 H^{N-i}(X,\C)_\chi((\hbar)),\qquad
 \bigoplus_a H^{N-2a-i}((X^\vee)^g_a,\C)^G((\hbar)).
\]
The second initially uses Azumaya trace lines with canonical finite
connections. Their lowest-degree cohomology on isolated charts forces
equivariant triviality (Lemma~\ref{qr:untwisting}). Since a fixed
component of age $a$ has dimension $N-2a$, Hochschild degree $i$
pairs source degree $k=N-i$ with sector degree $k-2a$.

To recover the Tate weights as well, we retain the cap action of
\[
 \kappa=\operatorname{KS}(t\partial_t)\in
 \HH^2(\cA_{\mathrm{gen}}).
\]
Under both realizations, a global Fedosov calculation identifies this
action with the derivative of the period
(Proposition~\ref{qr:cap-period}). Invariant second cohomology has
type $(2,2)$, so after a scalar rescaling this is cup product with
a formal deformation of the GSV class of the same type. Curious
Lefschetz identifies the resulting monodromy index $j=d-p$ with
source Tate degree $p$ and sector Tate degree $p-a$.

Writing $H^{k,p}(Z)=\Gr^W_{2p}H^k(Z,\C)$, we obtain
\[
 \dim H^{k,p}(X)_\chi
 =\sum_a\dim H^{k-2a,p-a}((X^\vee)^g_a)^G
\]
for every $k,p$ (Theorem~\ref{qr:character-comparison}). The split
Tate decompositions give the characterwise motivic comparison and its
rational sum. For odd $N$, adjoining an unused frozen coordinate
reduces to even dimension; its common $\Gm$ factor cancels in the
bigraded Tate multiplicity polynomials.

\subsection{Related work}

The form of the motivic statement parallels the fixed-sector comparison
of Hoskins and Pepin Lehalleur for Higgs bundles~\cite{HPL}. Our proof
uses quantum cluster algebras and symplectic deformation quantization
in place of the Hitchin fibration. The geometric and Hodge-theoretic
inputs are the cohomology of cluster varieties developed by Lam and
Speyer~\cite{LS}, the sink-recurrent extension in the work on braid
varieties~\cite{GLSBS,GLSB}, and Hlavinka's motivic splitting
theorem for RFR sink-recurrent cluster varieties over number
fields~\cite[Theorem~4.30]{Hlavinka}. Theorem~\ref{geom:splitting}
extends the required splitting to full-rank varieties and their
diagonal fixed components. The quantum localization argument uses
Muller's freezing results~\cite{MullerQuantum} and root specialization
as in Nguyen--Trampel--Yakimov~\cite{NTY}. The isolated calculation
starts from the generalized Weyl algebra resolution of
Solotar--Su\'arez-Alvarez--Vivas~\cite{SSAV}. The formal comparisons
and their operator compatibility use~\cite{Yekutieli,Dolgushev,DE,WX}.
The Azumaya-valued construction uses the deformation theory and
Fedosov resolutions of~\cite{BGNTAzumaya,DolgushevEndomorphisms}.
The classical trace lines are those of~\cite{TVdB}; the comparison
with the formal deformation and its Kodaira--Spencer action is
the content of Section~\ref{qr:section}.
The quotient presentation also places the construction alongside the
Cox-quotient approach to toric mirror symmetry in~\cite{ShendeToric}.

\subsection*{Acknowledgments}

We thank Mingyuan Hu for generously sharing the results of his joint
work with Vivek Shende and for helpful correspondence. The authors would also like to thank Josephine Hlavinka, Yifeng Huang, Marc Levine, Yeqin Liu and Peng Zhou 
for helpful discussions. This material is based upon work supported by the National Science Foundation 
under Grant No. DMS-2424139, while the first author was in residence at the Simons Laufer 
Mathematical Sciences Institute in Berkeley, California, during the Fall 2026 semester. Generative AI provided valuable feedback during the preparation of this manuscript and assisted in identifying and addressing gaps in the proofs. 
The authors take full responsibility for the correctness of all mathematical content.

\section{Cluster varieties and the mirror quotient}\label{geom:section}

Let $A=(B;C)\in M_{N\times n}(\Z)$, where $N=n+f$ and
$A$ has full column rank. Fix a positive integral diagonal matrix
$D=\operatorname{diag}(d_i)$ such that $DB=-B^{\mathsf T}D$,
taking $D=I_n$ when $B$ is skew-symmetric. Put
\begin{equation}\label{geom:dual-data}
 m=\operatorname{lcm}(d_i),\qquad p_i=m/d_i,\qquad
 P_m=\operatorname{diag}(p_i),\qquad P=\operatorname{diag}(P_m,I_f),
 \qquad A^\vee=\begin{pmatrix}-B^{\mathsf T}\\CP_m\end{pmatrix}.
\end{equation}
When there are no mutable indices, we take $m=1$.
The initial mutable variables are $x_1,\ldots,x_n$, and the frozen variables
$u_1,\ldots,u_f$ are invertible. The ordinary cluster algebra
$\cA(A)$ is generated inside $\C(x,u)$ by the frozen Laurent
ring and all cluster variables. The upper algebra $\cU(A)$ is
the intersection of the Laurent rings of all seeds. We write
$X(A)=\Spec\cA(A)$ and $X^\vee=X(A^\vee)$.
Our arrow convention is $b_{ij}>0$ for an
arrow $i\to j$. In direction $i$, the exchange ratio is
\[
 \widehat y_i=x^{B_{*,i}}u^{C_{*,i}}.
\]
The Laurent phenomenon~\cite{FZ} gives
\begin{equation}\label{geom:initial-torus}
 \cA(A)[(x_1\cdots x_n)^{-1}]
 =\C[x_1^{\pm1},\ldots,x_n^{\pm1},u_1^{\pm1},\ldots,u_f^{\pm1}].
\end{equation}

\subsection{Ordinary sink covers}

Freezing a mutable set deletes its exchange columns and makes its
variables invertible frozen variables; rows are reordered accordingly.

\begin{definition}\label{geom:sr}
The class of \emph{sink-recurrent} seeds is the smallest class
containing all seeds with zero principal matrix, closed under
mutation, and satisfying the following recursive condition:
a seed belongs to the class if it has a nonisolated sink $s$
with mutable neighbor set $T\ne\varnothing$ such that both
seeds obtained by freezing $\{s\}$ and $T$ belong to the class.
\end{definition}

The convention in \cite[Definition~3.6]{GLSB} freezes
$T\cup\{s\}$ in the second branch. It gives the same class:
after $T$ is frozen, $s$ is isolated, and adding or deleting an
isolated mutable vertex does not affect sink recurrence.
Thus the recursion gives a finite rooted tree with isolated leaves.

Full column rank is preserved by mutation and freezing. The same
holds for the really full-rank condition that
$A^{\mathsf T}:\Z^N\to\Z^n$ be surjective: mutation is multiplication
by integral unimodular matrices \cite[formulas~(3.1)--(3.3)]{BZ},
and freezing composes $A^{\mathsf T}$ with a coordinate projection.

The identity $PA^\vee=AP_m$ shows that $A^\vee$ has full rank.
Moreover, $-B^{\mathsf T}=DBD^{-1}$ has the same signs and zero
pattern as $B$, and this positive diagonal conjugacy commutes with
mutation. Hence $A$ and $A^\vee$ are sink-recurrent simultaneously.

\begin{lemma}[{\cite[Lemma~1]{MullerAU}}]\label{geom:freezing}
For the ordinary and upper algebras of a frozen seed $F$, one has
\[
 \cA_F\subseteq\cA(A)[x_i^{-1}:i\in F]\subseteq\cU_F.
\]
Thus $\cA_F=\cU_F$ implies that the frozen algebra is the
indicated ordinary localization.
\end{lemma}

\begin{proposition}\label{geom:geometry}
If $A$ is full-rank and sink-recurrent, then $\cA(A)=\cU(A)$,
every freezing used in the recursion of Definition~\ref{geom:sr} is the corresponding
ordinary localization, and $X(A)$ is smooth, integral, affine,
and of finite type, of dimension $N$.
\end{proposition}

\begin{proof}
For isolated seeds, ordinary--upper equality is
\cite[Proposition~3]{MullerAU}. At a recursive sink,
\eqref{geom:sink-cover} makes $x_s$ and $x_T$ comaximal.
Inductively, Lemma~\ref{geom:freezing} identifies the two frozen
algebras with ordinary localizations, and
\cite[Lemma~2]{MullerAU} gives ordinary--upper equality on their
union. Thus the recursion gives a finite cover by isolated cluster
localizations. In particular, $\cA(A)$ is locally acyclic.
Finite generation and smoothness now follow from
\cite[Theorems~4.2 and~7.7]{Muller}; these results also hold for
skew-symmetrizable matrices, as stipulated in \cite[\S1.1]{Muller}.
Integrality follows from $\cA(A)\subset\C(x,u)$, and
\eqref{geom:initial-torus} gives the dimension.
The sink-recurrent ordinary--upper equality is also recorded in
\cite[Proposition~3.7]{GLSB}.
\end{proof}

At an isolated seed, put $R=\C[u^{\pm1}]$ and let $c_i$ be the columns of $C$.
Rescaling exchange partners by frozen units gives
\begin{equation}\label{geom:isolated}
 \cA(A)=
 R[x_1,y_1,\ldots,x_n,y_n]/(x_iy_i-1-u^{c_i})_{i=1}^n.
\end{equation}
At a nonisolated sink $s$ with mutable neighbor set $T$, the normalized
exchange relation is
\begin{equation}\label{geom:sink-cover}
 x_s\widetilde x_s'=1+\eta\prod_{j\in T}x_j^{b_{js}},
 \qquad \eta\in\C[u^{\pm1}]^\times,\quad b_{js}>0.
\end{equation}
Thus $x_s$ and $x_T=\prod_{j\in T}x_j$ generate the unit ideal, giving
\[
 X(A)=D(x_s)\cup D(x_T).
\]
When $s$ is chosen as in Definition~\ref{geom:sr}, these opens are the
cluster varieties obtained by freezing $\{s\}$ and $T$, respectively.

\subsection{Coefficient changes and diagonal actions}

Coefficient rescaling \cite[Propositions~2.3 and~3.2]{Fraser} applies
to two geometric-type patterns with the same principal matrix and
frozen rings embedded in a common Laurent ring. If their initial
mutable variables differ by frozen monomials and their exchange ratios
agree, the same holds after every mutation word. Their ordinary
algebras over the common frozen ring therefore agree.
For rank one, the same assertion follows directly from the single
exchange relation.

The diagonalizable group
\begin{equation}\label{geom:grading-group}
 T_A=\ker\left(\Gm^N\longrightarrow\Gm^n,
        \ (a,b)\longmapsto
        (a^{B_{*,i}}b^{C_{*,i}})_{i=1}^n\right)
\end{equation}
acts on $X(A)$ by scaling its initial coordinates. Its character
group is $\Z^N/A\Z^n$. Write
\[
 S_A=\pi_0(T_A),\qquad
 \Xi_A=S_A^\vee=\sat(A\Z^n)/A\Z^n.
\]
In particular, $T_A$ is a connected torus of dimension $f$
if $A$ is really full-rank. This is the diagonal
action of \cite[Proposition~5.1]{LS}; its exchange-homogeneity
criterion applies unchanged to skew-symmetrizable matrices.

For the Hu--Shende mirror $\cM_A=[P_A/H]$ of~\eqref{intro:mirror}
\cite{HS}, the frozen torus of $P_A$ is
$T_0=\mathbb G_{m,Z}^n\times\mathbb G_{m,U}^f$. The action in
\eqref{intro:action} extends to all of $P_A$: its character
matrix satisfies
\[
 [I_n,B^{\mathsf T},C^{\mathsf T}]
 \begin{pmatrix}-B^{\mathsf T}\\I_n\\0\end{pmatrix}=0,
\]
so every exchange relation is homogeneous.

\subsection{The finite quotient presentation}

Put
\begin{equation}\label{geom:finite-matrices}
 S=B+P_mC^{\mathsf T}C,\qquad d=|\det S|,
 \qquad J_0=\begin{pmatrix}-B^{\mathsf T}&-C^{\mathsf T}\\CP_m&-I_f\end{pmatrix}.
\end{equation}
The matrix $S$ is nonsingular. Indeed, if $Sv=0$ for a real
vector $v$, then $0=v^{\mathsf T}DSv=m\|Cv\|^2$, and hence
$Cv=Bv=0$. Full rank of $A$ gives $v=0$. Also
$|\det J_0|=d$ by the Schur complement.

\begin{proposition}\label{geom:mirror}
Define
\begin{equation}\label{geom:original-group}
 G_0=\{(a,a^{-C})\in\Gm^{n+f}:
           a^{S_{i,*}}=1\text{ for all }i\}.
\end{equation}
Then $G_0\subset T_{A^\vee}$ is finite of order $d$, its character
group is $\Z^N/J_0\Z^N$, and there is an isomorphism of stacks
\begin{equation}\label{geom:mirror-quotient}
 \cM_A\simeq[X^\vee/G_0].
\end{equation}
\end{proposition}

\begin{proof}
Consider the isogeny of frozen tori
\begin{equation}\label{geom:frozen-isogeny}
 \begin{split}
 H_h\times\mathbb G_{m,u}^f&\longrightarrow T_0,\\
 Z_i&=h^{B_{i,*}}u^{p_iC_{*,i}},\qquad
 U_\ell=h^{C_{\ell,*}}u_\ell^{-1}.
 \end{split}
\end{equation}
Its exponent matrix is
$\left(\begin{smallmatrix}B&P_mC^{\mathsf T}\\C&-I_f\end{smallmatrix}\right)$,
whose determinant has absolute value $d$ by the Schur complement.
Thus it is finite \et\ of degree $d$.

We claim that its pullback to $P_A$ is
\begin{equation}\label{geom:mirror-base-change}
 P_A\times_{T_0}(H_h\times\mathbb G_{m,u}^f)
 \simeq H_h\times X^\vee.
\end{equation}
Let $R$ be the ordinary coordinate ring of $P_A$, let
$R_0=\C[Z^{\pm1},U^{\pm1}]$, and let
$R_1=\C[h^{\pm1},u^{\pm1}]$. The isogeny makes $R_1$
finite free over $R_0$, by representatives of the quotient
character lattice. Hence $R'=R\otimes_{R_0}R_1$ is finite free
over $R$, and multiplication by $X_1\cdots X_n$ on $R'$ is
injective. Its localization is the Laurent domain
$\C[X^{\pm1},h^{\pm1},u^{\pm1}]$. Thus $R'$ embeds in
that domain.

Set $X_i=h_ix_i$. The principal exchange ratio becomes
\[
 Z_i\prod_jX_j^{-b_{ij}}
 =h^{B_{i,*}}u^{p_iC_{*,i}}\prod_j(h_jx_j)^{-b_{ij}}
 =u^{p_iC_{*,i}}\prod_jx_j^{-b_{ij}}.
\]
Coefficient rescaling identifies $R'$ with
$\C[h^{\pm1}]\otimes\cA(A^\vee)$, proving
\eqref{geom:mirror-base-change}.
The resulting finite \et\ map is
\begin{equation}\label{geom:mirror-map}
 \Phi(h,x,u)=
 \bigl(X_i=h_ix_i,\ Z_i=h^{B_{i,*}}u^{p_iC_{*,i}},\
             U_\ell=h^{C_{\ell,*}}u_\ell^{-1}\bigr).
\end{equation}

The kernel of~\eqref{geom:frozen-isogeny} consists of
$(\varepsilon,\eta)$ with
$\eta_\ell=\varepsilon^{C_{\ell,*}}$ and
$\varepsilon^{S_{i,*}}=1$. On $H\times X^\vee$ its deck action
is
\[
 (h,x,u)\longmapsto
 (\varepsilon h,\varepsilon^{-1}x,\eta u).
\]
Putting $a=\varepsilon^{-1}$ identifies the action on the
second factor with~\eqref{geom:original-group}. Moreover
\[
 a^{-B_{i,*}}(a^{-C})^{p_iC_{*,i}}
 =a^{-(B+P_mC^{\mathsf T}C)_{i,*}}=1,
\]
so $G_0\subset T_{A^\vee}$. Its equations are equivalently those of
the diagonal kernel with character lattice $\Z^N/J_0\Z^N$.

The map $\Phi$ is $H$-equivariant, with $H$ multiplying $H_h$.
The $H$-action commutes with the deck group, so taking their
quotients in either order gives~\eqref{geom:mirror-quotient}.
\end{proof}

\subsection{Alternating completions and inertia motives}

\begin{lemma}\label{geom:alternating-completion}
If $N$ is even, there is an integral nonsingular alternating
matrix
\begin{equation}\label{geom:alternating-matrix}
 K=\begin{pmatrix}BP_m&-P_mC^{\mathsf T}\\CP_m&D_f\end{pmatrix}
   =PJ,\qquad
 J=\begin{pmatrix}-B^{\mathsf T}&-C^{\mathsf T}\\CP_m&D_f\end{pmatrix},
 \qquad D_f^{\mathsf T}=-D_f.
\end{equation}
More generally, for any full-column-rank matrix $F\in M_{N\times n}(\Z)$ and
nonsingular integral matrix $Q$ whose first columns are $F$,
the finite diagonal group
$G_Q=\operatorname{Hom}(\Z^N/Q\Z^N,\C^\times)$ is a subgroup
of $T_F$ and surjects onto $\pi_0(T_F)$.
\end{lemma}

\begin{proof}
Apply \cite[Lemma~5.12]{LS} to $AP_m$: its columns are independent
and its principal part $BP_m$ is skew-symmetric. The resulting
integral alternating completion has the required block form.

Put $L=F\Z^n\subset\Z^N$. Since the first columns of $Q$
are $F$, the character quotient defines $G_Q\subset T_F$.
Furthermore
\begin{equation}\label{geom:component-injection}
 \sat(L)\cap Q\Z^N=L.
\end{equation}
Indeed, if $Fr=Qz$ with $r\in\Q^n$ and $z\in\Z^N$, then
invertibility of $Q$ gives $z=(r,0)$, so $r\in\Z^n$.
Thus the torsion subgroup of $\Z^N/L$ injects into
$\Z^N/Q\Z^N$. Dualizing this injection proves the asserted
surjectivity onto the component group.
\end{proof}

The quantum comparison uses the auxiliary quotient $[X^\vee/G_J]$.
Both $G_J$ and $G_0$ surject onto $\pi_0(T_{A^\vee})$;
we identify their additive inertia motives.

\begin{lemma}\label{geom:torsion}
Let a connected complex algebraic torus $T$ act on a smooth
finite-type variety $V$. A finite-order point $g\in T(\C)$
acts as the identity on $M(V)$ with rational coefficients.
If the order of $g$ divides $e$, the same holds with coefficients
$\Z[1/e]$.
\end{lemma}

\begin{proof}
The unit point splits the torus motive as
$M(\Gm)=\Z\oplus\Z(1)[1]$. The reduced component of the
point morphism of $z\in\C^\times$ is its Kummer class in
\[
 \Hom(\Z,\Z(1)[1])=\C^\times;
\]
see \cite[Lemma~2.13 and Theorem~4.1]{MVW}. If $z^e=1$, this
class vanishes after inverting $e$. Split $T\simeq\Gm^r$.
Every coordinate of $g$ is such a root of unity, so the tensor
product formula gives $M(i_g)=M(i_1)$ for the point morphisms
to $T$. Composing with the motive of the action gives
\[
 M(g)=M(T\times V\to V)\circ(M(i_g)\otimes\id)
     =\id_{M(V)}.
\]
\end{proof}

\begin{proposition}\label{geom:inertia-comparison}
Assume that $A$ is sink-recurrent and $N$ is even, and choose
$K$ as in~\eqref{geom:alternating-matrix}. The potentially
nonempty fixed sectors for $G_0$ and $G_J$ have the same indices,
namely
\begin{equation}\label{geom:sector-group}
 H_A=\{(a,1):a^{B_{i,*}}=1,
                  \ a^{C_{\ell,*}}=1\text{ for all }i,\ell\}.
\end{equation}
For $h\in H_A$ and each age $a$, the averaging projectors for
$G_0$ and $G_J$ on $M((X^\vee)^h_a)$ are equal. Consequently
\begin{equation}\label{geom:inertia-identification}
 \Morb(\cM_A)\simeq
 \bigoplus_{h\in G_J}\ \bigoplus_a M((X^\vee)^h_a)^{G_J}\{a\}.
\end{equation}
\end{proposition}

\begin{proof}
An element with a fixed point must act trivially on every frozen
unit, and hence has the form $(a,1)$. For $G_0$, its equations
then reduce to~\eqref{geom:sector-group}. For $G_J$, the first
$n$ character equations give $a^{-B_{i,*}}=1$, and the remaining
ones give $a^{-C_{\ell,*}}=1$, which are the same conditions.
Thus the automorphisms with potentially nonempty fixed loci, and
hence their ages, agree.

The fixed loci are smooth because $X^\vee$ is smooth and the
automorphisms have finite order in characteristic zero. Each
union $V=(X^\vee)^h_a$ is $T_{A^\vee}$-stable. Let $G'$ be the finite
subgroup of $T_{A^\vee}$ generated by $G_0$ and $G_J$. Every element
of $G'\cap T_{A^\vee}^0$ is torsion, and acts trivially on $M(V)$
by Lemma~\ref{geom:torsion}. The action of $G'$ on $M(V)$
therefore factors through its image in $\pi_0(T_{A^\vee})$.
Both $G_0$ and $G_J$ surject onto that component group by
Lemma~\ref{geom:alternating-completion}. Their averaging
projectors are therefore the same component-group average.

The diagonal action preserves the symplectic GSV form associated
with $P^{-1}KP^{-1}$, whose extension to $X^\vee$ is established in
Section~\ref{sec:splitting}. Thus the ages are integers and
$2a$ is the codimension on each fixed component. Applying the
projector equality to these components and using
Proposition~\ref{geom:mirror} gives
\eqref{geom:inertia-identification}.
\end{proof}

\subsection{The geometric quotient in the primitive case}

\begin{proposition}\label{geom:slice}\label{geom:rfr}
Suppose that $A$ is really full-rank. Choose an integral matrix
$L$ such that
\[
 \mathcal V=\begin{pmatrix}A^{\mathsf T}\\L\end{pmatrix}
 \in\operatorname{GL}_N(\Z),\qquad
 C_{\mathrm{sl}}=L_{*,1:n},\qquad
 A_{\mathrm{sl}}=\begin{pmatrix}-B^{\mathsf T}\\C_{\mathrm{sl}}\end{pmatrix}.
\]
Then $A_{\mathrm{sl}}$ is really full-rank, and there is an
$H$-equivariant isomorphism
\[
 P_A\simeq H\times X(A_{\mathrm{sl}}),
\]
with $H$ acting by multiplication on the first factor. In
particular, $\cM_A\simeq X(A_{\mathrm{sl}})$ has trivial inertia.
If $B$ is skew-symmetric and $A$ is sink-recurrent, the quotient morphism
\begin{equation}\label{geom:rfr-map}
 p_A:X(A)\longrightarrow\cM_A\simeq X(A_{\mathrm{sl}}),
 \qquad (x,u)\longmapsto[X=x,\ Z=u^C,\ U=u^{-1}],
\end{equation}
is a finite \et\ $G_0$-torsor of degree $d$ between smooth
finite-type affine varieties.
\end{proposition}

\begin{proof}
Surjectivity of $A^{\mathsf T}$ permits the chosen unimodular
completion. On $T_0$, use coordinates $(h,v)\in H\times\Gm^f$
defined by
\[
 t_\ell=h^{A_{\ell,*}}v^{L_{*,\ell}},\qquad t=(Z,U).
\]
This is an isomorphism of tori with cocharacter matrix
$\mathcal V^{\mathsf T}$. Set $X_i=h_ix_i$. Then
\[
 Z_i\prod_jX_j^{-b_{ij}}
 =v^{(C_{\mathrm{sl}})_{*,i}}\prod_jx_j^{-b_{ij}}.
\]
Coefficient rescaling identifies the entire ordinary algebra
with $\C[h^{\pm1}]\otimes\cA(A_{\mathrm{sl}})$.
The $H$-action multiplies $h$ and fixes $x,v$, proving
equivariance. The first $n$ columns of $\mathcal V$ are
$(B^{\mathsf T};C_{\mathrm{sl}})$ and span a primitive direct summand.
Changing the signs of the first $n$ rows gives $A_{\mathrm{sl}}$.
Hence $A_{\mathrm{sl}}$ is really full-rank. Its principal matrix
$-B^{\mathsf T}$ is sink-recurrent whenever $B$ is.

When $B$ is skew-symmetric, our convention $D=I_n$ gives $X^\vee=X(A)$.
The free $H$-action and Proposition~\ref{geom:mirror} then show that
$G_0$ acts freely on $X(A)$ with quotient $X(A_{\mathrm{sl}})$.
Thus the quotient map is a finite \et\ torsor, with
formula~\eqref{geom:rfr-map} obtained by setting $h=1$ in
\eqref{geom:mirror-map}. Proposition~\ref{geom:geometry} gives
smoothness and finite type in the sink-recurrent case.
\end{proof}

\subsection{Comparison with the Fock--Goncharov ensemble}
\label{geom:fg}

To identify the Fock--Goncharov atlas inside $\cM_A$, choose the
initial square completion
\begin{equation}\label{geom:fg-seed}
 \widehat B=\begin{pmatrix}B&-P_mC^{\mathsf T}\\C&0\end{pmatrix},
 \qquad \varepsilon=\widehat B^{\mathsf T},\qquad
 \varepsilon^\vee=-\varepsilon^{\mathsf T}=-\widehat B,
\end{equation}
with the first $n$ indices mutable. In the Fock--Goncharov row
convention, $\varepsilon_{ka}=A_{ak}$, so its $\mathcal A$-atlas
has the exchange relations of $A$. Up to a common factor, the
multipliers for $\varepsilon$ and $\varepsilon^\vee$ are respectively
\[
 (p_1,\ldots,p_n,1,\ldots,1),\qquad
 (d_1,\ldots,d_n,m,\ldots,m).
\]
Let $\mathcal X_{\mathrm{FG}}(A)$ be the Langlands-dual
$\mathcal X$-variety obtained by gluing the seed tori along the
maximal isomorphism domains of composed mutations
\cite[\S1.2]{FG}, \cite[Definition~A.4]{GHKK}.
We mutate the full matrix~\eqref{geom:fg-seed}, including its
submatrix indexed by frozen variables.

\begin{proposition}\label{geom:fg-affinization}\label{geom:fg-mutation}
Mutation induces isomorphisms $\cM_A\simeq\cM_{\mu_k(A)}$
whose transitions on the quotient tori are the
Fock--Goncharov mutations for $\varepsilon^\vee$.
If $A$ is sink-recurrent, these tori define an open immersion
\[
 \mathcal X_{\mathrm{FG}}(A)\hookrightarrow\cM_A.
\]
If, in addition, $A$ is really full-rank, then
\begin{equation}\label{geom:fg-affine}
 \cM_A\simeq
 \Spec\Gamma(\mathcal X_{\mathrm{FG}}(A),\mathcal O).
\end{equation}
\end{proposition}

\begin{proof}
Write $A'=(B';C')=\mu_k(A)$.
On the initial cluster torus of $P_A$, the $H$-action is a
trivial torsor with quotient coordinates
\begin{equation}\label{geom:fg-coordinates}
 t_i=Z_i\prod_jX_j^{-b_{ij}},\qquad
 w_\ell=U_\ell\prod_jX_j^{-c_{\ell j}}.
\end{equation}
Mutate $X_k$ in $Y_{-B^{\mathsf T}}$, writing $T_i$ for the
resulting mutable variables, and change the frozen coordinates by
\begin{equation}\label{geom:fg-frozen-change}
 \begin{split}
 Z_k^{\mathrm{new}}&=Z_k^{-1},\\
 Z_j^{\mathrm{new}}&=Z_jZ_k^{[-b_{jk}]_+}\quad(j\ne k),\\
 U_\ell^{\mathrm{new}}&=U_\ell Z_k^{[-c_{\ell k}]_+}.
 \end{split}
\end{equation}
The mutated principal coefficient block $E$ is the identity
except in row $k$, where $E_{kk}=-1$ and
$E_{kj}=[-b_{jk}]_+$ for $j\ne k$. Since $E^2=I$, this Laurent change
restores the identity block. The exchange ratios are those of
$(-(B')^{\mathsf T};I_n;0)$, so coefficient
rescaling identifies the ordinary algebras.

The weights of the $T_i$ are the columns of the matrix $F$
which is the identity except in column $k$, where $F_{kk}=-1$
and $F_{ik}=[b_{ki}]_+$ for $i\ne k$. Put $a_\ell=[-c_{\ell k}]_+$.
Then $F^2=I$, and matrix mutation gives
\[
 F B^{\mathsf T}E=(B')^{\mathsf T},\qquad
 F(C^{\mathsf T}+B^{\mathsf T}e_ka^{\mathsf T})
       =(C')^{\mathsf T}.
\]
Changing the character basis of $H$ by $F$ therefore gives
weights $[I_n,(B')^{\mathsf T},(C')^{\mathsf T}]$, proving
that the algebra isomorphism descends to $\cM_A\simeq\cM_{A'}$.
The exchange relation
$X_kT_k=(\prod_jX_j^{[b_{kj}]_+})(1+t_k)$ now gives, for
$z=(t,w)$,
\begin{equation}\label{geom:fg-mutation-formula}
 z'_k=z_k^{-1},\qquad
 z'_a=z_a z_k^{[-A_{ak}]_+}(1+z_k)^{A_{ak}}\quad(a\ne k).
\end{equation}
Since $\varepsilon^\vee_{ak}=-A_{ak}$, these are precisely
\cite[formula~(13)]{FG}.

To identify the glued atlas, take principal coefficients for the
full dual $N$-element seed, with the same mutable indices. Set
\[
 Z=\Spec\cA\!\begin{pmatrix}-B^{\mathsf T}\\CP_m\\I_n\\0_{f\times n}\end{pmatrix},
\]
with initial coordinates $(x;v;\zeta;\omega)$ of lengths
$n,f,n,f$. The torus $T=\Gm^N$ acts by
\[
 q\cdot x_i=q_ix_i,\quad q\cdot v_\ell=q_{n+\ell}v_\ell,
 \quad q\cdot(\zeta,\omega)_a
       =q^{\widehat B_{a,*}}(\zeta,\omega)_a.
\]
The exchange ratios are invariant, so the action extends to
$Z$. Replacing $\zeta_i$ by
$\widehat\zeta_i=\zeta_i\prod_\ell v_\ell^{p_ic_{\ell i}}$
gives
\[
 Z\simeq\mathbb G_{m,v}^f\times Y_{-B^{\mathsf T}}
                     \times\mathbb G_{m,\omega}^f.
\]
The factor $\Gm^f\subset T$ acts only on $v$, while
$\Gm^n\subset T$ has weights $[I_n,B^{\mathsf T},C^{\mathsf T}]$
on $(x;\widehat\zeta;\omega)$. Thus
\begin{equation}\label{geom:fg-principal-quotient}
 [Z/T]\simeq\cM_A.
\end{equation}

The union $\mathcal A_{\mathrm{prin}}$ of the cluster tori is
$T$-stable and open in $Z$. Its torus intersections are the
maximal isomorphism domains of the composed mutations, since
an intersection is the locus where both sets of cluster
coordinates are invertible. By \cite[Proposition~B.2(2)]{GHKK},
$\mathcal A_{\mathrm{prin}}\to\mathcal X_{\mathrm{FG}}(A)$ is
a $T$-torsor. Quotienting the invariant open inclusion by $T$
and using~\eqref{geom:fg-principal-quotient} proves the open
immersion.

For sink-recurrent $A$, this principal pattern is full-rank and
sink-recurrent. Proposition~\ref{geom:geometry} gives
\[
 \Gamma(\mathcal A_{\mathrm{prin}},\mathcal O)=\mathcal O(Z).
\]
If $A$ is really full-rank, Proposition~\ref{geom:slice} makes
$\cM_A$ affine, and descent of functions gives
\[
 \Gamma(\mathcal X_{\mathrm{FG}}(A),\mathcal O)
 =\Gamma(\mathcal A_{\mathrm{prin}},\mathcal O)^T
 =\mathcal O(Z)^T=\mathcal O(\cM_A).
\]
This proves~\eqref{geom:fg-affine}.
\end{proof}

\begin{remark}
For skew-symmetric $B$ and really full-rank $A$, the map $p_A$
of~\eqref{geom:rfr-map} and the ensemble map $p_{\mathrm{ens}}$
satisfy, on the initial source torus,
\cite[formula~(11)]{FG}
\[
 \begin{aligned}
 p_A^*t_i&=p_{\mathrm{ens}}^*t_i=x^{B_{*,i}}u^{C_{*,i}},\\
 p_A^*w_\ell&=u_\ell^{-1}p_{\mathrm{ens}}^*w_\ell
              =u_\ell^{-1}x^{-C_{\ell,*}}.
 \end{aligned}
\]
The factors $u_\ell^{-1}$ come from the finite isogeny defining $p_A$.
\end{remark}

\section{Split Tate motives and curious Lefschetz}
\label{sec:splitting}

We use the notation $H^{k,p}=\Gr^W_{2p}H^k$. A mixed Tate
variety of dimension $2d$ satisfies \emph{curious Lefschetz} for a
class $\omega\in H^{2,2}$ if
\[
 \omega^{d-p}:H^{k,p}\xrightarrow{\sim}
 H^{k+2(d-p),\,2d-p}\qquad(p\leq d).
\]
\begin{theorem}\label{geom:splitting}
Let $X=X(A)$ be full-rank and sink-recurrent. Its motive is a finite
split Tate motive:
\[
 M(X)\simeq\bigoplus_{k,p}\Q(p)[k]^{\oplus h^{k,p}(X)},
 \qquad h^{k,p}(X)=\dim_\Q H^{k,p}(X).
\]
Every connected component of a fixed locus of a finite-order diagonal
cluster automorphism has the same property. These assertions hold over
a sufficiently large cyclotomic number field before base change to
$\C$, and also for invariant summands under finite diagonal groups.
\end{theorem}

For really full-rank sink-recurrent cluster varieties, splitting
over number fields is due to Hlavinka
\cite[Theorem~4.30]{Hlavinka}. We extend the argument to full-rank
seeds and fixed components, using curious Lefschetz and arithmetic splitting.

Suppose $N$ is even. A rational alternating matrix $K$ is
\emph{compatible} with $A$ if
$Ke_i=q_iA_{*,i}$ for nonzero rational numbers $q_i$ at the mutable
indices. Let $K$ be such a matrix and suppose it is nonsingular.
The form
\begin{equation}\label{geom:gsv}
 \omega_K=\sum_{a<b}K_{ab}\,d\log x_a\wedge d\log x_b
\end{equation}
is a compatible form in the sense of Gekhtman--Shapiro--Vainshtein
\cite[Theorem~2.1]{GSV2005}. Its mutation invariance is the GSV
transformation rule; see also
\cite[Proposition~6.3 and Corollary~6.5]{FG}. The rule applies
linearly to the compatible coefficients. Compatibility is also
preserved by freezing.

The isolated charts cover $X$. Independence of their frozen exchange
columns gives a finite cover of the frozen torus on which the
equations become a product of
\begin{equation}\label{geom:rankone}
 X_m=\{xy=1+z^m,\ z\ne0\}
\end{equation}
and a torus, as in~\cite[Proposition~5.10]{LS}. The complement of the cluster
tori has codimension at least two, since the omitted locus in each
surface is $x=y=0$; compare \cite[Proposition~7.6]{LS}.
Thus~\eqref{geom:gsv} extends to a regular closed
form on the smooth variety $X$. Its top exterior power is nonzero on
all cluster tori. A zero of a nonzero section of the canonical line
bundle would have a divisorial component; hence the form is symplectic
everywhere. As a rational logarithmic form, its class has type
$(2,2)$~\cite[Proposition~2.13]{LS}.

For $K$ in~\eqref{geom:alternating-matrix}, the matrices $K$ and
$P^{-1}KP^{-1}$ are compatible on $X$ and
$X^\vee$, with weights $p_i$ and $p_i^{-1}$, respectively.

\begin{lemma}\label{geom:isolated-lefschetz}
At a full-rank isolated seed, every nondegenerate rational compatible
form satisfies curious Lefschetz. Its class, divided by $(2\pi i)^2$,
is rational.
\end{lemma}
\begin{proof}
After a finite frozen cover the exchange matrix is $(0;dI_n;0)$.
In the corresponding coordinates the form matrix is
\[
 \begin{pmatrix}0&-L&0\\ L&S&-R^{\mathsf T}\\0&R&T\end{pmatrix},
 \qquad L=d\operatorname{diag}(q_i).
\]
Apply the frozen-monomial normalization in the proof of
\cite[Proposition~8.2]{LS}. The same triangular changes work with
the invertible diagonal matrix $L$ in place of $dI_n$; after a
further finite power cover they are integral and eliminate $S$ and $R$.
Nondegeneracy makes $T$ nonsingular. The form is then a sum
of nonzero rational multiples of $d\log x_i\wedge d\log z_i$ on the
surfaces~\eqref{geom:rankone} and a nondegenerate rational alternating
form on the remaining torus. Rationality and curious Lefschetz follow
from the rank-one, torus, and product calculations in
\cite[Propositions~7.1, 8.1, and 3.4]{LS}. Finite pullback is injective
on rational cohomology and identifies the original cohomology with
the deck-invariant part. The invariant form therefore descends both
assertions.
\end{proof}

\begin{lemma}\label{geom:frozen-h1}
The classes $(2\pi i)^{-1}d\log u_\ell$ of the frozen coordinates
form a basis of $H^1(X,\Q)$.
\end{lemma}
\begin{proof}
Restriction to the initial torus is injective in degree one.
A class from $H^1(X)$ has zero residue along every boundary divisor
visible on an $i$-mutated torus. The one-exchange residue calculation
in the proof of~\cite[Proposition~9.4]{LS} gives residue $c_i$
for the representative $\sum_a c_a d\log x_a$; full rank makes
the exchange character nonconstant, so these divisors are nonempty.
Thus $c_i=0$ for every mutable index. The frozen classes extend
globally and are independent on the initial torus.
\end{proof}

At a sink, the three opens $U=D(x_s)$, $V=D(x_T)$, and $U\cap V$
have the same dimension as $X$, and the restrictions of $\omega_K$
are compatible and nondegenerate. Curious Lefschetz propagates from
them to $X$ by~\cite[Theorem~3.5]{LS}. Induction from
Lemma~\ref{geom:isolated-lefschetz} is the sink-cover argument
of~\cite[Theorem~10.1]{GLSBS}; the skew-symmetrizable extension is
also noted in \cite[\S6.4, footnote~6]{GLSB}. Moreover,
\[
 H^1(U,\Q)\oplus H^1(V,\Q)\longrightarrow H^1(U\cap V,\Q)
\]
is surjective by Lemma~\ref{geom:frozen-h1}: the frozen coordinates
on the intersection are the union of those on the two opens.
Thus $H^2(X,\Q)$ injects into $H^2(U,\Q)\oplus H^2(V,\Q)$.
Its complexification sends $(2\pi i)^{-2}[\omega_K]$ to a rational
class by induction. Since the injection is defined over $\Q$, its
unique preimage is rational as well.
Curious Lefschetz for this rational class implies that all cohomology
is split mixed Tate~\cite[Theorem~3.3]{LS}.

\begin{proposition}\label{geom:fixed-lefschetz}
Let $N=2d$, and let $h$ be a finite-order diagonal cluster
automorphism. Every connected component $F$ of $X^h$ satisfies curious
Lefschetz for $\omega_K|_F$, and has split mixed Tate cohomology.
Its dimension is even and
\[
 \age(h,F)=\tfrac12(N-\dim F).
\]
\end{proposition}
\begin{proof}
If the fixed locus is nonempty, $h$ fixes every frozen coordinate.
It preserves $\omega_K$. A finite-order symplectic automorphism has
smooth symplectic fixed loci; its nontrivial normal eigenvalues occur
in inverse pairs, proving the age formula.

Every cluster coordinate is an eigenfunction, so the sink tree consists
of $h$-stable opens. The smooth connected component $F$ is irreducible.
Its intersection with a nonempty open in that tree is therefore a
connected component of the fixed locus there. These intersections
are affine: they are principal opens in the closed affine subvariety
$F$. Restrict the tree to $F$, omitting empty branches.
At an isolated leaf write $x_ix_i'=1+u^{C_i}$ and let $J$ index the
nontrivial mutable eigenvalues. The fixed equations are
\[
 x_j=x_j'=0,\qquad u^{C_j}=-1\quad(j\in J).
\]
On a component of this frozen torus fiber put $u=u_0v^R$, where $u_0$
is torsion and the columns of $R$ are an integral basis of
$\ker(C_J^{\mathsf T})$. The remaining exchange columns are
$C'=R^{\mathsf T}C_{J^c}$. They are independent, since a dependence
would put a combination of the remaining columns of $C$ in the span
of $C_J$. The associated morphism of tori is surjective, so a further
torsion translation removes the root-of-unity constants in the
remaining exchange relations. This component is a full-rank isolated
cluster variety.

The restriction of the form has matrix
\[
 \begin{pmatrix}0&-Q'(C')^{\mathsf T}\\ C'Q'&R^{\mathsf T}D_fR\end{pmatrix},
 \quad
 K=\begin{pmatrix}0&-QC^{\mathsf T}\\CQ&D_f\end{pmatrix},\qquad
 Q=\operatorname{diag}(q_i),\quad Q'=\operatorname{diag}(q_i)_{i\in J^c}.
\]
It is the nondegenerate symplectic reduction by the isotropic space
spanned by the removed mutable directions. It computes the actual
restriction across their vanishing locus, because the regular identity
\[
 d\log x_j\wedge d\log u^{C_j}
   =\frac{dx_j\wedge dx_j'}{x_jx_j'-1}
\]
shows, after multiplication by $q_j$, that each removed summand
restricts to zero. Lemma
\ref{geom:isolated-lefschetz} gives curious Lefschetz at the fixed
leaves. The open-cover theorem propagates it up the restricted tree;
all its nonempty opens have dimension $\dim F$. The already established
global rational class restricts rationally, so~\cite[Theorem~3.3]{LS}
also gives Hodge splitting on $F$.
\end{proof}

\begin{lemma}\label{geom:arithmetic-splitting}
Let $k$ be a number field and $M$ a geometric mixed Tate motive over
$k$. If the cohomology of every Hodge realization of $M$, for every
embedding $k\hookrightarrow\C$, is split mixed Tate, then $M$ is a
finite direct sum of shifted Tate motives.
\end{lemma}
\begin{proof}
The category of mixed Tate motives over a number field has its bounded
Tate $t$-structure; its heart has homological dimension one and its
bounded derived category is the geometric mixed Tate category. The
family of realization maps
\[
 \operatorname{Ext}^1_{\mathrm{MT}(k)}(\Q(a),\Q(b))\longrightarrow
 \bigoplus_{\sigma:k\hookrightarrow\C}
 \operatorname{Ext}^1_{\mathrm{MHS}_{\Q}}(\Q(a),\Q(b))
\]
is injective. These number-field statements and the all-embeddings
regulator are recorded in~\cite[\S3.4, (87)--(88), and \S3.5,
Lemma~3.4, (92)]{Goncharov}. As in
\cite[Proposition~4.27]{Hlavinka}, induct on the weight filtration
of a cohomology object. Once the lower-weight part is split, the
next extension is a matrix of Tate extensions. Its image at every
embedding vanishes because the Hodge realization is split and its
weight filtration is canonical.
Injectivity of the displayed map makes the extension itself zero.
Homological dimension one splits the derived object into its shifted
cohomology objects.
Tensoring the extension groups and realization maps with a number
field $E$ proves the same assertion with $E$ coefficients.
\end{proof}

\begin{proof}[Proof of Theorem~\ref{geom:splitting}]
Choose the finite tree of sink covers, including the intersections in
each cover. On every isolated chart, stratify by the vanishing mutable
coordinates. For a set $I$ of vanishing mutable coordinates, the
stratum in~\eqref{geom:isolated} is
\[
 \A^{|I|}\times\mathbb G_m^{n-|I|}
       \times\{u\in\mathbb G_m^f:u^{c_i}=-1\ (i\in I)\}.
\]
The affine coordinates are the exchange partners $y_i$, for $i\in I$;
the other exchange partners are determined. Independence of the
columns makes the last factor a fiber of a surjective torus morphism.
Over one finite cyclotomic extension $k$, the finitely many connected
strata have the form
\[
 S_{\lambda,k}\simeq\A_k^{a_\lambda}\times\mathbb G_{m,k}^{b_\lambda}.
\]
Localization and smooth duality give mixed Tate motives for the
isolated charts. The Mayer--Vietoris triangles
\[
 M(U\cap V)\longrightarrow M(U)\oplus M(V)
 \longrightarrow M(U\cup V)\longrightarrow M(U\cap V)[1]
\]
then show that $M(X_k)$ is geometric mixed Tate.

Since $X$ is defined over $\Q$, every conjugate
$X_\sigma=X_k\times_{k,\sigma}\C$ is $X$. Thus
\[
 H^j(X_\sigma,\Q)\simeq
 \bigoplus_p\Q(-p)^{\oplus h^{j,p}(X)}
 \qquad\text{for every }\sigma.
\]
Lemma~\ref{geom:arithmetic-splitting} therefore splits $M(X_k)$.

Now let $F$ be a connected component of $X^h$, with $h$ of finite
order $r$. Enlarge $k$ to define $h$ and the connected components on
the isolated charts. The closure of each such component in $X^h$ is a
global component, so $F$ is defined over $k$. Restricting the same
stratifications and covers to $F$ shows that $M(F_k)$ is mixed Tate.
For every embedding $\sigma:k\hookrightarrow\C$,
\[
 \sigma(h)=h^{a_\sigma},\qquad
 \gcd(a_\sigma,r)=1,\qquad X^{h^{a_\sigma}}=X^h.
\]
Thus the conjugates of $F$ are components of the same fixed locus.
Proposition~\ref{geom:fixed-lefschetz} gives Hodge splitting for every
conjugate, and Lemma~\ref{geom:arithmetic-splitting} splits $M(F_k)$.

Finite-group invariant summands satisfy the same hypotheses: taking
invariants is a rational idempotent operation, so it preserves both
mixed Tate motives and split Hodge structures. If $N$ is odd, apply
the even-dimensional argument after adjoining one unused frozen
coordinate. The decomposition
\[
 M(X\times\Gm)\simeq M(X)\oplus M(X)(1)[1]
\]
and its analogue for the fixed loci exhibit the original motives as
retracts. Their Hodge realizations are split at every embedding, so
the arithmetic lemma applies to these retracts as well. Finally,
Hodge realization determines the Tate multiplicities by degree and
weight, giving the displayed formula in the theorem.
\end{proof}

\section{Quantum localization and completion}\label{ql:section}

We prove Hochschild descent for sink covers over the parameter ring
and after completion. This reduces the comparison at two formal
specializations to isolated seeds.

Let $E$ be a field of characteristic zero and $\widetilde B$ an
extended exchange matrix of full column rank with skew-symmetrizable
principal part. Fix an integral alternating matrix $\Lambda$ satisfying
\[
 \widetilde B^{\mathsf t}\Lambda=(\Delta\ \ 0),
\]
where $\Delta$ is a diagonal matrix with positive integer entries. All frozen
variables are inverted. Our based quantum-torus convention is
\begin{equation}\label{qm:torus}
 M_t(v)M_t(w)=t^{\Lambda(v,w)}M_t(v+w).
\end{equation}
Write
$\mathcal A_R$ and $\mathcal U_R$ for the ordinary and upper quantum
cluster algebras over $R$, respectively. Here $R$ is
$E[t^{\pm1}]$, its fraction field, or a localization of
$E[t^{\pm1}]$ at a nonzero prime. The ordinary algebra is generated
by the quantum cluster variables and the frozen inverses; the upper
algebra is the intersection of the quantum tori over all seeds.
Both are regarded as subalgebras of the same division algebra.
We use the quantum Laurent phenomenon of
\cite[Corollary~5.2]{BZ}; for a strict specialized seed we use
\cite[Theorem~3.10]{NTY}.

\subsection{Sink localization}

We use the recursion in Definition~\ref{geom:sr}, with the two
freezings $\{s\}$ and $T$ at each nonisolated sink $s$, where
$T$ is its mutable neighbor set.

\begin{lemma}\label{ql:comaximal}
Let $A$ be a ring, and suppose that $x,y\in A$ satisfy
$xy=c\,yx$ for a central unit $c$ and $Ax+Ay=A$. Then
\begin{equation}\label{ql:powers}
 Ax^r+Ay^s=A\qquad(r,s\geq1).
\end{equation}
Suppose, in addition, that the algebras obtained by adjoining
$x^{-1}$, $y^{-1}$, and both inverses are two-sided Ore
localizations $S_x,S_y,S_{xy}$ of $A$ inside a common division
algebra. Then there is an exact sequence of $A$-bimodules
\begin{equation}\label{ql:cech-bimodule}
 0\longrightarrow A
 \xrightarrow{\ a\mapsto(a,a)\ } S_x\oplus S_y
 \xrightarrow{\ (u,v)\mapsto u-v\ } S_{xy}
 \longrightarrow0.
\end{equation}
\end{lemma}

\begin{proof}
Consider the left $A$-module $A/(Ax^r+Ay^s)$ and let $v$ be the
class of $1$. Quasi-commutation implies that $x^iy^jv=0$ if
$i\geq r$ or $j\geq s$. If $v\ne0$, choose a nonzero such
monomial of maximal total degree, and denote it by $w$. Then
$xw=yw=0$. An expression $1=ax+by$ gives $w=0$, a contradiction.

If $z\in S_x\cap S_y$, choose $r,s$ such that $x^rz,y^sz\in A$.
Writing $1=ax^r+by^s$ shows that $z\in A$. For $z\in S_{xy}$,
choose $r,s$ such that $x^ry^sz\in A$. Then $x^rz\in S_y$ and
$y^sz\in S_x$, again by quasi-commutation, and the same expression
for $1$ shows that $z\in S_x+S_y$. These assertions prove
exactness of the displayed bimodule sequence.
\end{proof}

\begin{proposition}\label{ql:ordinary-upper}
Suppose that the principal part of $\widetilde B$ is
sink-recurrent. Then $\mathcal A_R=\mathcal U_R$. Every freezing
used in this recursion is an actual two-sided Ore
localization of the ordinary quantum cluster algebra.
The same assertions hold for the strict specialized seed at any
$t=\alpha\in E^\times$.

More precisely, at a nonisolated sink $s$ occurring in this recursion, put
\[
 x=X_s,\qquad
 y=\prod_{j\in T}X_j^{|b_{js}|},
\]
with any fixed ordering of the product. The three algebras obtained
by freezing $s$, $T$, and $T\cup\{s\}$ are $S_x,S_y,S_{xy}$,
respectively, and satisfy~\eqref{ql:cech-bimodule} with
$A=\mathcal A_R$.
\end{proposition}

\begin{proof}
We use Muller's isolated-seed and freezing results
\cite[Propositions~8.2, 8.5, and~8.6]{MullerQuantum} and induct along
Definition~\ref{geom:sr}. The proofs carry over to compatible
skew-symmetrizable seeds because they use the exchange relations and
quasi-commutation. They also
apply over the stated coefficient rings and at $t=\alpha$, since the
only inverted coefficients are powers of $t$ and frozen monomials,
and the exchange binomials remain nonzero in the specialized quantum
tori. The isolated case is therefore covered by
\cite[Proposition~8.2]{MullerQuantum}.

For the inductive step, mutate to a seed with a nonisolated sink
$s$ for which both prescribed freezings are sink-recurrent, and
write $A=\mathcal A_R$. Since $s$ is a sink, one mutable part of its
exchange relation is $1$, and the other is $y$. The frozen factors
are units, so the exchange relation gives
\[
 Ax+Ay=A.
\]
Also $xy=c\,yx$ for a power $c$ of $t$.

For a set $F$ of mutable vertices, denote its frozen algebra by
$\mathcal A_R^{(F)}$. The freezing inclusions
\cite[Proposition~8.5]{MullerQuantum} are
\begin{equation}\label{ql:freezing-sandwich}
 \mathcal A_R^{(F)}
 \subseteq A[X_j^{-1}:j\in F]
 \subseteq\mathcal U_R[X_j^{-1}:j\in F]
 \subseteq\mathcal U_R^{(F)}.
\end{equation}
The middle expressions denote generated subalgebras of the ambient
division algebra. Induction identifies the two outer algebras when
$F=\{s\}$ or $F=T$. Hence the corresponding generated subalgebras
are precisely the ordinary frozen algebras. The same holds for
$F=T\cup\{s\}$: after $T$ is frozen, $s$ is isolated, and deleting
an isolated mutable vertex preserves sink recurrence.

These equalities give Ore localizations by
\cite[Proposition~8.6]{MullerQuantum}; the fraction argument with
right denominators also applies after specialization. Since the variables indexed
by $T$ pairwise quasi-commute, inverting $y$ is equivalent to
inverting every one of them.

Lemma~\ref{ql:comaximal} now proves~\eqref{ql:cech-bimodule} and
$S_x\cap S_y=A$. Every element of $\mathcal U_R$ belongs to
both frozen upper algebras, which induction has identified with
$S_x$ and $S_y$. Thus $\mathcal U_R\subseteq A$, completing the induction.
\end{proof}

\subsection{Hochschild descent}

For an $R$-flat algebra $A$, let
\[
 C_\bullet(A/R)=A\otimes^{\mathbf L}_{A\otimes_R A^{\mathrm{op}}}A
\]
be its relative Hochschild complex. We regard this as an object of
the derived category $\mathrm D(R)$, with
$\HH_i(A/R)=H^{-i}C_\bullet(A/R)$.

\begin{proposition}\label{ql:hochschild-descent}
At every sink step of Proposition~\ref{ql:ordinary-upper}, the
localization maps induce a distinguished triangle
\begin{equation}\label{ql:hochschild-triangle}
 C_\bullet(A/R)\longrightarrow
 C_\bullet(S_x/R)\oplus C_\bullet(S_y/R)
 \longrightarrow C_\bullet(S_{xy}/R)
 \longrightarrow C_\bullet(A/R)[1].
\end{equation}
The same assertion holds after reduction modulo any nonzero ideal
of a coefficient discrete valuation ring $R$.
\end{proposition}

\begin{proof}
A two-sided flat Ore localization $A\to S$ is a homological ring
epimorphism: multiplication induces
$S\otimes_A^{\mathbf L}S\simeq S$. Put
$A^e=A\otimes_R A^{\mathrm{op}}$ and similarly for $S^e$.
Extension of scalars gives
\[
 S^e\otimes_{A^e}^{\mathbf L}A\simeq S,
 \qquad
 S\otimes_{A^e}^{\mathbf L}A
 \simeq S\otimes_{S^e}^{\mathbf L}S.
\]
Applying $(-)\otimes_{A^e}^{\mathbf L}A$ to
\eqref{ql:cech-bimodule} proves~\eqref{ql:hochschild-triangle}.
The arrows are the usual localization maps on Hochschild complexes.

All terms of~\eqref{ql:cech-bimodule} embed in quantum Laurent
algebras, so they are torsion-free and hence $R$-flat.
Reduction modulo an ideal therefore preserves exactness and
two-sided flat ring epimorphisms. The same argument over the
quotient coefficient ring proves the last assertion.
\end{proof}

\subsection{Specialization and the formal algebra}

For the rest of this section, assume that $\widetilde B$ is
sink-recurrent. Fix a root of unity $\alpha\in E^\times$, and put
\[
 R_\alpha=E[t^{\pm1}]_{(t-\alpha)},\qquad
 \pi=t-\alpha.
\]
Write $\mathcal A_\alpha$ for the strict root-of-unity quantum
cluster algebra obtained by specializing the given compatible
integral seed at $t=\alpha$. In particular,
$\mathcal A_1$ is the ordinary commutative cluster algebra.

\begin{proposition}\label{ql:specialization}
The specialization homomorphism induces an isomorphism
\begin{equation}\label{ql:special-fiber}
 \mathcal A_{R_\alpha}/\pi\mathcal A_{R_\alpha}
 \xrightarrow{\ \sim\ }\mathcal A_\alpha.
\end{equation}
The $\pi$-adic completion of $\mathcal A_{R_\alpha}$ is
topologically free over $E[[\pi]]$, with special fiber
$\mathcal A_\alpha$. These assertions also hold for every frozen
algebra appearing in this recursion.
\end{proposition}

\begin{proof}
We apply the proof of the specialization theorem
\cite[Theorem~5.11]{NTY} over $R_\alpha$, using
Proposition~\ref{ql:ordinary-upper}. Its divisibility input
\cite[Proposition~5.10]{NTY} applies over $R_\alpha$ and gives, for
adjacent quantum tori $\mathcal T$ and $\mathcal T'$,
\begin{equation}\label{ql:adjacent-divisibility}
 \mathcal T\cap\pi\mathcal T'
 =\pi\mathcal T\cap\mathcal T'.
\end{equation}
The proof there requires only that the exchange binomials remain
nonzero modulo $\pi$ in the common quantum subtorus, which holds
because the exchange columns are nonzero. Thus the saturation
conclusion used in \cite[Theorem~5.11]{NTY} is
\begin{equation}\label{ql:saturation}
 \mathcal A_{R_\alpha}\cap\pi\mathcal T_0
 =\pi\mathcal A_{R_\alpha},
\end{equation}
where $\mathcal T_0$ is the initial quantum torus. The integral
compatible matrix $\Lambda$ supplies the strict lift, and
\cite[Theorem~5.7]{NTY} identifies the image with
$\mathcal A_\alpha$.

Put $M=\mathcal A_{R_\alpha}$ and $V=M/\pi M$. Since $M$ is flat
over $R_\alpha$, the quotients $M/\pi^rM$ form an adic system of
flat $E[[\pi]]/(\pi^r)$-modules. By
\cite[Proposition~2.4]{Yekutieli}, their inverse limit is
topologically free with special fiber $V$:
\begin{equation}\label{ql:topologically-free}
 V[[\pi]]\xrightarrow{\ \sim\ }\widehat M.
\end{equation}
Here $\widehat M=\varprojlim_rM/\pi^rM$, and the isomorphism is
noncanonical. The same argument applies to the frozen algebras.
\end{proof}

Use the formal parameter $\hbar$ defined by
\begin{equation}\label{ql:formal-parameter}
 t=\alpha\exp(\hbar).
\end{equation}
Since $\pi=\alpha\hbar+O(\hbar^2)$, this identifies
$E[[\pi]]$ with $E[[\hbar]]$. Denote the completed algebra by
$\widehat{\mathcal A}_\alpha$. It is a complete, flat,
topologically free associative deformation of
$\mathcal A_\alpha$. The completion is taken in the parameter:
modulo each power of $\hbar$, an element is an ordinary algebra
element with finite support.

\subsection{Continuous Hochschild chains}

Write
\[
 \mathcal A_r=\mathcal A_{R_\alpha}/\pi^r\mathcal A_{R_\alpha}.
\]
We identify $R_\alpha/\pi^r$ with
$E[[\hbar]]/(\hbar^r)$ using~\eqref{ql:formal-parameter}.
Define
\begin{equation}\label{ql:continuous-definition}
 C_\bullet^{\mathrm{cont}}(\widehat{\mathcal A}_\alpha)
 =\mathbf R\!\varprojlim_r
 C_\bullet(\mathcal A_r/(R_\alpha/\pi^r)).
\end{equation}
The derived inverse limit is taken in $\mathrm D(E[[\hbar]])$.
This agrees with the completed bar construction: each finite-level
algebra is flat over its coefficient ring, and the transition maps
between the bar-chain modules are surjective. Write
\[
 C_\bullet^{\mathrm{cont}}(\widehat{\mathcal A}_\alpha)[\hbar^{-1}]
 =C_\bullet^{\mathrm{cont}}(\widehat{\mathcal A}_\alpha)
   \otimes_{E[[\hbar]]}E((\hbar)).
\]

\begin{proposition}\label{ql:continuous-descent}
Completing the algebras at a sink step gives a distinguished
triangle
\begin{align}\label{ql:continuous-triangle}
 C_\bullet^{\mathrm{cont}}(\widehat A_\alpha)
 &\longrightarrow
 C_\bullet^{\mathrm{cont}}(\widehat S_{x,\alpha})
 \oplus C_\bullet^{\mathrm{cont}}(\widehat S_{y,\alpha})
 \longrightarrow
 C_\bullet^{\mathrm{cont}}(\widehat S_{xy,\alpha})
 \longrightarrow
 C_\bullet^{\mathrm{cont}}(\widehat A_\alpha)[1].
\end{align}
It remains a distinguished triangle after inverting $\hbar$.
\end{proposition}

\begin{proof}
Proposition~\ref{ql:hochschild-descent} gives a compatible
triangle at every finite level. Derived inverse limits preserve
finite distinguished triangles. Extension of scalars to
$E((\hbar))$ preserves them as well.
\end{proof}

Let $\mathcal A_\eta=\mathcal A_{E(t)}$ and use
\eqref{ql:formal-parameter} to embed $E(t)$ in $E((\hbar))$.
The map from algebraic bar chains to completed bar chains induces
a canonical comparison
\begin{equation}\label{ql:comparison-map}
 \gamma_\alpha\colon
 C_\bullet(\mathcal A_\eta/E(t))
 \otimes_{E(t)}E((\hbar))
 \longrightarrow
 C_\bullet^{\mathrm{cont}}(\widehat{\mathcal A}_\alpha)
 [\hbar^{-1}].
\end{equation}
It commutes with the cluster localizations above.

\begin{theorem}\label{ql:completion}
Suppose that $\widetilde B$ is full-column-rank and
sink-recurrent and that $\Lambda$ is nondegenerate. Then
$\gamma_\alpha$ is a quasi-isomorphism. In particular, for every
$i$ it induces an isomorphism
\begin{equation}\label{ql:completion-homology}
 \HH_i(\mathcal A_\eta/E(t))
 \otimes_{E(t)}E((\hbar))
 \xrightarrow{\ \sim\ }
 H^{-i}\!\left(
 C_\bullet^{\mathrm{cont}}(\widehat{\mathcal A}_\alpha)
 [\hbar^{-1}]\right).
\end{equation}
\end{theorem}

\begin{proof}
For an isolated seed this is Theorem~\ref{qt:completion}.
Its bounded chain contractions establish precisely the
canonical map~\eqref{ql:comparison-map}. We induct on the number
of mutable vertices. Full column rank persists under mutation and
freezing. Mutation changes $\Lambda$ by an integral congruence
\cite[Proposition~3.4]{BZ},
and freezing retains the same alternating form on the full seed
lattice, so nondegeneracy also persists.

At a nonisolated sink, apply
Proposition~\ref{ql:hochschild-descent} over $E(t)$ and extend
scalars to $E((\hbar))$. Naturality of the bar-chain map gives a
morphism from this triangle to the triangle of
Proposition~\ref{ql:continuous-descent}, after inverting
$\hbar$. The two frozen algebras and their common freezing have
fewer mutable vertices and are sink-recurrent. Their comparison
maps are quasi-isomorphisms by induction. The comparison at the
original algebra is therefore a quasi-isomorphism as well.
\end{proof}

For a finite diagonal grading action, the denominators are homogeneous
and all comparison maps are equivariant. Thus
Theorem~\ref{ql:completion} holds on every character summand over
a splitting field.

\section{The root fiber and the dual cluster variety}
\label{sec:quantum-twist}

Assume that $A$ is sink-recurrent and $N$ is even. Retain the notation of
\eqref{geom:dual-data} and the quantum family of
Section~\ref{ql:section}. The alternating completion $K$ has first
$n$ columns $A\operatorname{diag}(p_i)$, and
\[
 K=PJ,\qquad
 J=\begin{pmatrix}-B^{\mathsf T}&-C^{\mathsf T}\\
 C\operatorname{diag}(p_i)&D_f\end{pmatrix},\qquad
 \Lambda=-\ell K^{-1}.
\]
Here $D_f$ is alternating, $\Lambda$ is integral, and
$A^{\mathsf T}\Lambda=(\ell\operatorname{diag}(p_i)^{-1},0)$.
Choose a positive multiple $s$ of $2\ell$, and put
\begin{equation}\label{qm:root-data}
 T=sP,\qquad \rho=\exp(\pi i/(\ell s)),\qquad
 Q=T\Z^N/sK\Z^N\simeq\Z^N/J\Z^N,\qquad G=Q^\vee.
\end{equation}
All coordinate algebras in this section are taken over a cyclotomic
number field $E$ containing $\rho$ and the characters of $Q$;
base change to $\C$ gives the varieties of Section~\ref{geom:section}.
The group $G$ acts on $X^\vee$ through the auxiliary quotient associated with $J$.

The root commutator on the source torus is
\begin{equation}\label{qm:pairing}
 \beta_\rho(v,w)=\rho^{2\Lambda(v,w)}
       =\exp(-2\pi i\,v^{\mathsf T}K^{-1}w/s).
\end{equation}
Its radical is $sK\Z^N$. The choice of $s$ gives the stronger
identity
\begin{equation}\label{qm:isotropic}
 \rho^{\Lambda(Ta,Tb)}=1\qquad(a,b\in\Z^N).
\end{equation}
For $\bar w\in Q$, define a parameter-independent action by
\begin{equation}\label{qm:source-action}
 \theta_{\bar w}(M_t(v))=\beta_\rho(Tw,v)M_t(v).
\end{equation}
This preserves exchange relations, since
$\beta_\rho(Tw,A_i)=1$. Set
\[
 \mathcal B_t=\cA_t\mathbin{\#}Q.
\]
Write $u_{\bar w}$ for the group basis, so that
$u_{\bar w}a=\theta_{\bar w}(a)u_{\bar w}$.
The dual action of $G$ fixes $\cA_t$ and sends
$u_{\bar w}$ to $\gamma(\bar w)u_{\bar w}$ for $\gamma\in G$.

\subsection{The central dual cluster algebra}

Adapting the root-of-unity binomial calculation of
\cite[Proposition~4.4 and Theorem~4.6]{NTY}, we construct the dual
central algebra using symmetrizer-dependent powers.

\begin{proposition}\label{qm:frobenius}
In the initial seed, the monomials
\begin{equation}\label{qm:central-frame}
 z_a=M_\rho(Ta)u_{\bar a}^{-1}\qquad(a\in\Z^N)
\end{equation}
form a commutative central toric frame in the crossed product
of the initial quantum torus with $Q$.
Under mutation these frames satisfy the ordinary exchange
relations for $A^\vee$. They define an embedding
$\mathcal O(X^\vee)\hookrightarrow Z(\mathcal B_\rho)$,
compatible with the cluster localizations.
\end{proposition}

\begin{proof}
Equation~\eqref{qm:isotropic} and the crossed-product relation
give $z_az_b=z_{a+b}$ and centrality. At a mutable index $i$,
write the mutated variable as $Y+Z$, where
\[
 Y=M_\rho(-e_i+[A_i]_+),\qquad
 Z=M_\rho(-e_i+[-A_i]_+).
\]
Compatibility gives $\Lambda(A_i,-e_i+[-A_i]_+)=-\ell/p_i$,
and hence $YZ=\rho^{-2\ell/p_i}ZY$.
Since $\rho^{2\ell/p_i}=\exp(2\pi i/(sp_i))$ has order exactly
$sp_i$, the quantum binomial formula gives
$(Y+Z)^{sp_i}=Y^{sp_i}+Z^{sp_i}$.

If $E_i$ and $E_i^\vee$ are the source and dual mutation
involutions, then
\[
 E_iT=TE_i^\vee,\qquad
 sp_iA_i=T A_i^\vee=sK_i.
\]
The two surviving monomials therefore have exactly the dual
exchange exponents; their group factors agree, since
$A_i^\vee$ is a column of $J$. The mutated completion satisfies
$K'=E_iK E_i^{\mathsf T}$ and
$J'=E_i^\vee J E_i^{\mathsf T}$.
Thus $E_i^\vee$ induces an isomorphism
$\Z^N/J'\Z^N\longrightarrow\Z^N/J\Z^N$, under which the
group factor of the mutated central monomial of exponent $a$ is
$u_{\overline{E_i^\vee a}}^{-1}$.
Equation~\eqref{qm:isotropic} continues to hold for every
mutated integral form $\Lambda'$, because $T=sP$ and
$s\in2\ell\Z$. Thus the same calculation applies in every
seed. The map is injective on the initial Laurent torus, and
the Laurent phenomenon proves injectivity on the ordinary
cluster algebra. Inverting a dual coordinate in
\eqref{qm:central-frame} is equivalent to inverting its source
coordinate, which proves compatibility with localization.
\end{proof}

\begin{remark}
Liu identifies the parameters $v$ and $-v$ using a quadratic
refinement modulo two, odd on mutable directions
\cite[Lemmas~3.12--3.13 and Proposition~5.3]{LiuCoulomb}.
For~\eqref{qm:central-frame}, the divisibility $s\in2\ell\Z$
instead gives $\rho^{\Lambda'(Ta,Tb)}=1$ in every seed.
\end{remark}

\begin{proposition}\label{qm:center}
The algebra $\mathcal B_\rho$ is finite and reflexive over its
center, and
\[
 Z(\mathcal B_\rho)=\mathcal O(X^\vee).
\]
\end{proposition}

\begin{proof}
First work on a seed torus. The radical of
\eqref{qm:pairing} is $sK\Z^N$. Commutation with the source
monomials shows that a central term $M_\rho(v)u_{\bar a}$
has $v+Ta\in sK\Z^N$. By \eqref{qm:isotropic}, it then
commutes with $Q$ as well. Hence the center is precisely the
Laurent algebra of the frame~\eqref{qm:central-frame}.
The crossed product is finite free over this center.
Every geometric fiber is a twisted group algebra of a finite
group: its generators are invertible, and their powers are
prescribed nonzero scalars. The same homogeneous commutator
calculation gives its center as the ground field. Maschke's
theorem makes this algebra semisimple, so it is a matrix algebra.
The torus order is consequently Azumaya.

By Proposition~\ref{ql:ordinary-upper}, the strict source
algebra is the intersection of its seed tori. Taking a crossed
product with the fixed finite group $Q$ commutes with this
intersection, coefficient by coefficient in the group basis.
The corresponding ordinary--upper equality for $A^\vee$
therefore identifies the global center with
$\mathcal O(X^\vee)$.

Let $U\subset X^\vee$ be the union of the initial torus and
its single mutations. Full rank implies the starfish condition
\cite[Corollary~1.9]{BFZ}. Every codimension-one point lies in one
of these charts by \cite[Corollary~3.7]{Pompili}, so the complement
of $U$ has codimension at least two.
The torus orders glue to an Azumaya algebra
$\mathscr B_U$ on $U$. As $X^\vee$ is smooth, the direct image
\[
 \mathscr B=j_*\mathscr B_U,\qquad j:U\hookrightarrow X^\vee,
\]
is a coherent reflexive sheaf of algebras
\cite[Tag~0EBJ]{StacksReflexive}.
On any other seed torus, this sheaf and its torus order are
reflexive and agree away from a subset of codimension at least
two. They consequently agree on that torus. It follows that
\[
 \Gamma(X^\vee,\mathscr B)
   =\bigcap_{\text{seeds}}(\mathcal T_\rho\mathbin{\#}Q)
   =\mathcal B_\rho.
\]
Since $X^\vee$ is affine, coherence and reflexivity of
$\mathscr B$ prove the assertion.
\end{proof}

\subsection{The Azumaya algebra and finite duality}

\begin{proposition}\label{qm:azumaya}
The algebra $\mathcal B_\rho$ is Azumaya over $X^\vee$.
Its formal deformation at $t=\rho e^\hbar$ induces the
nondegenerate Poisson tensor $2T\Lambda T$ on the central
cluster torus.
\end{proposition}

\begin{proof}
Proposition~\ref{ql:specialization} gives a flat formal
deformation of $\mathcal B_\rho$ by taking the crossed product
with $Q$. For a central element $z$, choose a lift $\widetilde z$
linearly and set
\[
 D_z(b)=\hbar^{-1}[\widetilde z,\widetilde b]\bmod\hbar.
\]
These are derivations of $\mathcal B_\rho$, and their
restrictions give a Poisson bracket on its center. Thus
$\mathcal B_\rho$ is a Poisson order in the sense of
\cite[Sections~2.1--2.2]{BrownGordon}; different linear lifts
change the derivations by inner derivations.
On the central frame the bracket is
\[
 \{z_a,z_b\}=2a^{\mathsf T}T\Lambda Tb\,z_{a+b}.
\]
Its inverse is the weighted GSV form on $X^\vee$, which is
nondegenerate everywhere by the argument following
\eqref{geom:gsv}.
The smooth connected variety $X^\vee$ therefore has a single
symplectic leaf, hence a single symplectic core.

Over $\C$, Theorem~4.2 of \cite{BrownGordon} identifies the fibers
of $\mathcal B_\rho$ at all closed points as algebras. On the
initial torus they are matrix algebras of the same size.
The fiber-rank function of the coherent module $\mathscr B$ is
upper semicontinuous, and every nonempty closed subset of
$X^\vee_{\C}$ contains a closed point. Its rank therefore equals
the generic rank everywhere. Thus $\mathscr B$ is locally free
by~\cite[Tag~0FWG]{StacksReflexive}. The map
\[
 \mathscr B\otimes\mathscr B^{\mathrm{op}}
       \longrightarrow\mathcal E nd_{\mathcal O_{X^\vee}}(\mathscr B)
\]
is an isomorphism on every fiber, hence globally. The resulting
Azumaya property descends from $\C$ to $E$ by faithful flatness.
\end{proof}

\begin{theorem}\label{qm:morita}
The root crossed product $\mathcal B_\rho$ is an Azumaya
algebra on $X^\vee$, with a $G$-equivariant formal deformation
at $t=\rho e^\hbar$. For the whole parameter family there is
a parameter-independent isomorphism
\begin{equation}\label{qm:double-crossed}
 (\cA_t\mathbin{\#}Q)\mathbin{\#}G
       \simeq\operatorname{Mat}_{|Q|}(\cA_t).
\end{equation}
It induces a Morita equivalence compatible with the
Kodaira--Spencer class and its Hochschild cap action.
In the crossed-product decomposition of Hochschild homology,
the $\gamma$-sector on the left corresponds to the
$\gamma$-eigenspace for the source action of $Q$.
\end{theorem}

\begin{proof}
The Azumaya assertion is Proposition~\ref{qm:azumaya}, and
the dual $G$-action is defined over the parameter ring.
Cohen--Montgomery duality \cite[Theorem~3.2]{CohenMontgomery}, in the form of
\cite[Theorem~2.2]{NastasescuRodino}, identifies the left side of
\eqref{qm:double-crossed} with
\[
 \operatorname{End}_{\cA_t}
       ((\cA_t\mathbin{\#}Q)_{\cA_t}).
\]
Here the Fourier isomorphism $E[G]\simeq E[Q]^*$ identifies
the second crossed product with the smash product in that theorem.
The right module is free on the group basis, and the duality maps
are defined over the parameter ring using the fixed $Q$-action.
Thus the isomorphism is parameter-independent in these coordinates.
Differentiating identifies the Kodaira--Spencer cocycles; the matrix
Morita maps preserve cap products over $E(t)$
\cite[Theorem~2.2]{ArmentaKeller}. Their finite formulas also apply
to continuous chains after completion.

The sector decomposition also preserves this cap action. Work over
$E(t)$ or a formal Laurent field, and write $\mathcal B$ for the
corresponding realization of $\mathcal B_t=\cA_t\#Q$.
Write $v_\gamma$ for the $G$-basis and choose a $G$-invariant
parameter connection on $\mathcal B$.
Fixing the group basis extends its defect $c_{\mathcal B}$ to
\[
 c_{\#}(b v_\gamma,b'v_\eta)
   =c_{\mathcal B}(b,\gamma(b'))v_{\gamma\eta}.
\]
The bar sector maps of~\cite{Lorenz} move the group elements
to the coefficient factor. Substitution in this formula identifies
the crossed-product cap action with that of $c_{\mathcal B}$ on
each $\mathcal B\gamma$-twisted complex. The finite group sums
and averaging commute with completion.

For the character assertion, extend $\theta_q$ to the double crossed product
by fixing both group bases. Conjugation by $u_q$ acts on
$\cA_t$ as $\theta_q$ and sends $v_\gamma$ to
$\gamma(q)^{-1}v_\gamma$. Consequently
$\theta_q=\operatorname{Ad}(u_q)\circ\eta_q$, where
$\eta_q$ fixes $\cA_t\mathbin{\#}Q$ and multiplies
$v_\gamma$ by $\gamma(q)$. Inner automorphisms act trivially
on Hochschild homology, so the action on that sector is
precisely $\gamma(q)$.
\end{proof}

\subsection{Characters and fixed sectors}

The action of $\gamma\in G$ on the center is
\[
 \gamma(z_a)=\gamma(\bar a)^{-1}z_a.
\]
Thus its diagonal coordinates on $X^\vee$ are the inverses
of its values on the standard basis of $\Z^N/J\Z^N$.
With this convention, its image is exactly $G_J$.

\begin{lemma}\label{qm:finite-pairing}
The restriction of characters gives an injection
\[
 \Xi_A\hookrightarrow G.
\]
For $v=Ar\in\Z^N$, the corresponding character of $Q$ and
diagonal element of $G_J$ are respectively
\begin{equation}\label{qm:sector-label}
 \chi_r(\bar w)=\exp(-2\pi i\,r\cdot w_{\mathrm{mut}}),
 \qquad
 h_r=(e^{2\pi ir_1},\ldots,e^{2\pi ir_n},1_f).
\end{equation}
The image consists of the elements fixing all frozen
coordinates. In particular, every element with a nonempty
fixed locus has this form. The action of $Q$ on $X(A)$
surjects onto the component group of its diagonal grading torus.
\end{lemma}

\begin{proof}
The character represented by $v=Ar$ restricts along
\eqref{qm:source-action} as
\[
 \beta_\rho(Tw,v)
       =\exp(-2\pi i\,r\cdot w_{\mathrm{mut}}).
\]
It is trivial for every $w$ exactly when $r\in\Z^n$.
This proves injectivity, and finite diagonalizable duality
gives surjectivity onto the component group.
An element of $G_J$ with frozen coordinates one can be written
$h=(e^{2\pi ir},1_f)$. The equations $h^J=1$ say precisely
that $Br$ and $Cr$ are integral, or $Ar\in\Z^N$.
This proves the image description. Frozen coordinates are
units, so every element with a fixed point must fix them.
\end{proof}

On every sink-cover chart, retain the same groups and restrict their
actions. Inverting a central coordinate in~\eqref{qm:central-frame}
inverts its source coordinate, so~\eqref{qm:double-crossed} and
the sector labels restrict to these charts.
The trace local systems of the Azumaya algebra are shown to
be equivariantly trivial in Lemma~\ref{qr:untwisting}.

\section{Hochschild homology and the Tate grading}\label{qr:section}

The two formal realizations of the generic quantum cluster algebra
carry more information than their Hochschild gradings.  The
Kodaira--Spencer class of the quantum parameter acts on both
realizations.  We identify this action with multiplication by a
deformation of the GSV class and recover the Tate grading from its
monodromy filtration.

Throughout this section, $N=2d$, $X=X(A)$, and $X^\vee=X(A^\vee)$.
The finite source group is $Q=T\Z^N/sK\Z^N$, and
$G=Q^\vee=G_J$ acts diagonally on $X^\vee$.
Retain the splitting number field $E\subset\C$ of
Section~\ref{sec:quantum-twist}.
We work over $\C$ when using de Rham cohomology.  Write
$H^{k,(p,p)}(X)$ for the Deligne summand of the split mixed Tate
structure on $H^k(X,\C)$. The symplectic forms $\omega$ on $X$ and
$\omega^\vee$ on $X^\vee$ are inverse to the first-order Poisson
tensors of their quantizations. In logarithmic torus coordinates
their matrices are
\[
 -\frac{K}{2\ell},\qquad
 -\frac{P^{-1}KP^{-1}}{2\ell s^2},
\]
respectively. These are the corresponding weighted GSV forms.
All completed Hochschild complexes below use completed tensor products.

\subsection{The invariant part of the second cohomology}

\begin{lemma}\label{qr:invariant-second}
Let $A$ have full column rank, and suppose that $X(A)$ is smooth,
its ordinary and upper cluster algebras agree, and its cohomology
is mixed Tate.  Let $G$ be a finite subgroup of the diagonal cluster
automorphism group that surjects onto its group of connected
components.  Then
\[
 \operatorname{CH}^1(X(A))_{\Q}^{G}=0,
 \qquad
 H^2(X(A),\C)^G=H^{2,(2,2)}(X(A))^G.
\]
\end{lemma}

\begin{proof}
Put $M=\Z^N$ and $L=A\Z^n$.  For each exchange column $A_i$, let
$d_i$ be the greatest common divisor of its entries and put
$u_i=A_i/d_i$.  The class of $u_i$ in $M/L$ has order exactly $d_i$.
Indeed, if $e u_i=\sum_j a_j A_j$ with $a_j\in\Z$, column
independence over $\Q$ gives $a_i=e/d_i$ and $a_j=0$ for $j\ne i$.
The hypothesis on $G$ therefore implies that the character $u_i$
maps $G$ onto $\mu_{d_i}$.

Full rank implies the starfish condition.  The class-group
description in \cite[Theorems~3.3--3.4, Corollary~3.7 and
Theorem~3.8]{Pompili} consequently gives the following presentation
by boundary divisors:
\begin{equation}\label{qr:divisor-presentation}
 \operatorname{CH}^1(X)_{\Q}
 =\frac{\displaystyle\bigoplus_{i,\,\zeta^{d_i}=-1}
          \Q[D_{i,\zeta}]}{
        \displaystyle\left\langle
          \sum_{\zeta^{d_i}=-1}[D_{i,\zeta}]:i\text{ mutable}
        \right\rangle}.
\end{equation}
Here $D_{i,\zeta}$ is the boundary divisor corresponding to the
factor $x^{u_i}-\zeta$ of the exchange polynomial. Smoothness
identifies the divisor class group in the cited presentation with
$\operatorname{CH}^1(X)$.

For fixed $i$, the group $G$ acts transitively on the divisors
$D_{i,\zeta}$.  Taking invariants in
\eqref{qr:divisor-presentation} is exact over $\Q$, and each
invariant orbit sum is a displayed principal-divisor relation.
Thus $\operatorname{CH}^1(X)_{\Q}^{G}=0$.
The cycle map surjects onto the weight-two part of $H^2(X,\Q)$:
this follows from a smooth compactification and the Lefschetz
$(1,1)$ theorem.  Averaging gives the same assertion on
$G$-invariants.  Finally, the weight range for smooth varieties
and the mixed Tate hypothesis leave only weights two and four
in $H^2(X)$, proving the assertion.
\end{proof}

The lemma applies to the action of $Q$ on $X$ and to that of $G$
on $X^\vee$. Surjectivity onto the component groups follows from
Lemma~\ref{qm:finite-pairing} and
Lemma~\ref{geom:alternating-completion}, respectively. The other
hypotheses follow from the cluster localizations and
Theorem~\ref{geom:splitting}. In particular,
\begin{equation}\label{qr:second-purity}
 H^2(X)^Q=H^{2,(2,2)}(X)^Q,\qquad
 H^2(X^\vee)^G=H^{2,(2,2)}(X^\vee)^G.
\end{equation}

\subsection{The Kodaira--Spencer action under completion}

Let $\mathcal A$ be the generic quantum cluster algebra over $E(t)$
and put $\delta=t\partial_t$.  A $\delta$-connection on its
underlying vector space is an $E$-linear map $\nabla$ satisfying
$\nabla(fa)=\delta(f)a+f\nabla(a)$.  Its multiplication defect
\[
 c_\nabla(a,b)
 =\nabla(ab)-\nabla(a)b-a\nabla(b)
\]
is an $E(t)$-bilinear Hochschild $2$-cocycle.  Its class
$\kappa\in\HH^2(\mathcal A)$ is independent of $\nabla$.
We use the resulting cap action
$\kappa\cap-:\HH_i(\mathcal A)\to\HH_{i-2}(\mathcal A)$.

\begin{lemma}\label{qr:ks-completion}
For either embedding $t=e^{\hbar}$ or $t=\rho e^{\hbar}$,
the canonical comparison of Theorem~\ref{ql:completion}
intertwines cap product with $\kappa$ and the continuous
Kodaira--Spencer class of the formal family.
This compatibility is unchanged by an $\hbar$-dependent
equivalence of formal deformations.
\end{lemma}

\begin{proof}
Under either embedding, $\delta$ extends to $\partial_\hbar$.
Denote the map from the generic algebra to its formal realization
by $j$.  Choose a $\delta$-connection $\nabla$ on the generic
algebra and a continuous $\partial_\hbar$-connection
$\widehat\nabla$ on the formal algebra.  The latter exists by
a topological trivialization of the formal deformation; its
continuity permits a decrease of one in the $\hbar$-adic order.
On the generic algebra the difference
\[
 T=\widehat\nabla\circ j-j\circ\nabla
\]
is $E(t)$-linear with values in the formal algebra.  If $b$ denotes
the Hochschild coboundary with these coefficients, direct
expansion gives
\begin{equation}\label{qr:connection-defect}
 c_{\widehat\nabla}(j(a),j(b))
       -j\bigl(c_\nabla(a,b)\bigr)=-(bT)(a,b).
\end{equation}
The Hochschild cap identity makes the two induced maps from
algebraic bar chains to continuous bar chains chain homotopic.
The canonical comparison therefore intertwines their actions on
homology.  Formula~\eqref{qr:connection-defect}, applied to an
isomorphism of formal algebras, proves the assertion for an
$\hbar$-dependent equivalence as well.
\end{proof}

\subsection{Formal realizations and trace local systems}

At $t=e^\hbar$, the completed algebra is a $Q$-equivariant
quantization of $\mathcal O(X)$. At $t=\rho e^\hbar$,
Theorem~\ref{qm:morita} identifies the completed crossed product
$\mathcal A_t\#Q$ with a $G$-equivariant deformation
$\mathscr B_\hbar$ of an Azumaya algebra on $X^\vee$. The
parameter-independent matrix equivalence
\[
 (\mathcal A_t\#Q)\#G
 \simeq\operatorname{Mat}_{|Q|}(\mathcal A_t)
\]
recovers the original family and its Kodaira--Spencer action.

We give a common algebraic model for the two realizations. Let
$Y$ be smooth affine and symplectic, let $\Pi$ be a finite abelian
symplectic group, and let $\mathscr D_\hbar$ be a flat formal
deformation of a $\Pi$-equivariant Azumaya algebra $\mathscr D$ on
$Y$. Its first-order bracket on the center is assumed to be the
inverse $\pi$ of $\omega_Y$. All equivalences of deformations in
what follows induce the identity modulo $\hbar$.

The cochains of $\Gamma(Y,\mathscr D)$ that are multidifferential
over $Y$ include quasi-isomorphically into all Hochschild cochains.
Indeed, the Azumaya diagonal is perfect over its
enveloping algebra: \et\ matrix splittings reduce this to the
smooth diagonal of $Y$. Its derived endomorphisms therefore localize
along $Y$. Differential cochains, in turn, are filtered unions of
quasi-coherent sheaves. On matrix charts the inclusion induces the
HKR identification by cotrace and
\cite[Corollary~2.9 and Theorem~4.5]{YekutieliContinuous}; affine
acyclicity gives the global comparison. Its cohomology is
$\Gamma(Y,\bigwedge^\bullet T_Y)$ with the Schouten bracket.
This supplies the affine algebraic version of the Azumaya comparison
in~\cite[Theorem~34 and Corollary~35]{BGNTAzumaya}.
Choose a unit-preserving topological trivialization and average it
over $\Pi$. Finite invariants preserve the cochain quasi-isomorphism,
so the Maurer--Cartan comparison
\cite[Section~6 and Theorem~8.2]{Yekutieli} makes the actual product
differential while retaining its classical $\Pi$-action.

\begin{lemma}\label{qr:fedosov-model}
The actual deformation $\mathscr D_\hbar$ is equivariantly
equivalent to an Azumaya-valued algebraic Fedosov quantization.
It has a model on the Weyl bundle
\[
 \mathcal W_{\mathscr D}
 =\widehat{\operatorname{Sym}}(T_Y^*)[[\hbar]]
       \otimes_{\mathcal O_Y}\mathscr D
\]
with differential
\[
 D=\nabla-\delta+\hbar^{-1}\operatorname{ad}(r),
 \qquad \delta^{-1}r=0,
\]
and central curvature $-\Omega_\hbar/\hbar$, where
\[
 \Omega_\hbar=\omega_Y+\sum_{j\geq1}\hbar^j\Omega_j,
 \qquad \Omega_j\in\Gamma(Y,\Omega_Y^2)^{\Pi},
 \qquad d\Omega_j=0.
\]
The action on this model is the classical action on $\mathscr D$
and on the tangent bundle.
\end{lemma}

\begin{proof}
Choose an invariant torsion-free algebraic symplectic connection on $T_Y$
and an invariant algebra connection $\nabla^{\mathscr D}$ on
$\mathscr D$. They exist because their affine spaces of local
connections are torsors under quasi-coherent sheaves; these
torsors have no first cohomology on $Y$. Averaging preserves
the required properties. Every central-linear derivation of
$\mathscr D$ is inner. Consequently
\[
 (\nabla^{\mathscr D})^2=\operatorname{ad}(F_{\mathscr D})
\]
for a unique reduced-trace-zero two-form $F_{\mathscr D}$.
The Bianchi identity and reduced trace imply
$\nabla^{\mathscr D}F_{\mathscr D}=0$.
The product connection on $\mathcal W_{\mathscr D}$ has curvature
$\hbar^{-1}\operatorname{ad}(R_W+\hbar F_{\mathscr D})$.
The Fedosov recursion with contracting operator $\delta^{-1}$
therefore applies with $R_W+\hbar F_{\mathscr D}$ in place of
$R_W$. Give fiber coordinates degree one and $\hbar$ degree two.
The recursion increases this degree and produces $r$ of degree
at least three. It also preserves the property that $r$ modulo
$\hbar$ is scalar-matrix-valued: the noncentral initial
curvature is $\hbar F_{\mathscr D}$. Hence
$[r,\mathcal W_{\mathscr D}]\subseteq\hbar\mathcal W_{\mathscr D}$,
and $\hbar^{-1}\operatorname{ad}(r)$ preserves the integral
$\hbar$-adic lattice. Closedness of $\Omega_\hbar$ and the Bianchi
identity give $D^2=0$. Evaluation at the zero section identifies
the flat sections with $\mathscr D[[\hbar]]$ and their product
is differential. This is the connection construction underlying
\cite[Sections~3--4]{BGNTAzumaya}; for the simultaneous chain
and cochain resolution in the split case see
\cite[Section~3]{DolgushevEndomorphisms}.

We check that this construction includes the given deformation.
Start with $\Omega_\hbar=\omega_Y$. Its first-order Hochschild
class is $\pi$, so a differential gauge matches the products
modulo $\hbar^2$. Suppose that they agree modulo $\hbar^m$,
where $m\geq2$, and let $\nu_m$ be their difference in degree
$m$. Associativity in degrees $m$ and $m+1$ gives, respectively,
\[
 b\nu_m=0,\qquad [\pi,[\nu_m]_{\mathrm{HKR}}]=0.
\]
For the second equality one takes the Hochschild class of the
degree $m+1$ associativity equation; its remaining term is a
Hochschild coboundary. Since $\pi$ is symplectic,
$[\nu_m]_{\mathrm{HKR}}$ corresponds to a closed algebraic
two-form. Linearizing the Fedosov recursion shows that changing
$\Omega_\hbar$ by $\hbar^{m-1}\beta$ changes the degree $m$
Hochschild class by $-\pi\beta\pi$, with the fixed
antisymmetrization convention. In degree $m$ its HKR symbol
comes from the first fiber jets and one Moyal contraction;
each additional contraction supplies another factor of
$\hbar$. The symbol is therefore the linear term of
$(\omega_Y+\hbar^{m-1}\beta)^{-1}$. Choose $\beta$ to
match $[\nu_m]_{\mathrm{HKR}}$. The remaining difference is
$bP_m$ for a normalized differential one-cochain, and a gauge
$1+\hbar^mP_m$ removes it. The two-form and the one-cochain
may both be averaged over $\Pi$. These gauges converge
$\hbar$-adically, proving equivalence with the actual product,
with its actual finite action.
\end{proof}

For $g\in\Pi$ of order $q$, put
\[
 \mathscr L_g=
 \bigl(Z(\mathscr D|_{Y^g}\mathbin{\#}\langle g\rangle)\bigr)_g.
\]
This is the centralizer-equivariant line with its canonical
finite connection of~\cite[Corollary~1.29(i) and
Proposition~8.1]{TVdB}. Its multiplication gives
$\mathscr L_g^{\otimes q}\simeq\mathcal O_{Y^g}$.
Write $\mathscr L_{g,\hbar}=\mathscr L_g\otimes_\C\C((\hbar))$
for the corresponding local system after extension of scalars.
The fixed-sector comparison is
\begin{equation}\label{qr:twisted-hh}
 \HH_i^{\mathrm{cont}}(\mathscr D_\hbar,
                              \mathscr D_\hbar g)
 \simeq
 \bigoplus_{F\in\pi_0(Y^g)}
 H^{\dim F-i}(F,\mathscr L_{g,\hbar}),
\end{equation}
where global sections of the algebra are understood, $\hbar$ is
inverted, and the right bimodule action is twisted by $g$.
For $g=1$ this is the constant local system. We now prove the
comparison and the identification of its connection.

Let $F\subset Y^g$ be a connected component and put $d_F=\dim F$.
Take derived global sections of the normalized completed fiberwise bar complex
of the Fedosov model restricted to $F$, with total degree $p-j$
in exterior degree $p$ and chain degree $j$. Its differential
is $D+b$, with the total-complex signs. The direct sum of these
complexes over $F$ receives the flat-lift comparison from the
completed $g$-twisted bar complex. Modulo $\hbar$,
the Koszul differential $(1-g)y$ contracts the moving normal
directions, while the remaining jets give the Grothendieck
resolution of $\Omega_F^j\otimes\mathscr L_g$.
The comparison is twisted Azumaya HKR on \et\ matrix charts.
Its cone is derived $\hbar$-complete, so its vanishing modulo
$\hbar$ proves the comparison before and after inverting $\hbar$.
These are the jet comparisons of~\cite[Propositions~10--11]{DE}
with the matrix trace factor. All formal jet and parameter
completions here are derived inverse limits. The Weyl product
uses the combined filtration $\deg y=1$, $\deg\hbar=2$;
pure jet truncation is used only modulo $\hbar$.

After inverting $\hbar$, the relative fiber homology vanishes
outside chain degree $d_F$ by the calculation
of~\cite[Proposition~8]{DE} and matrix separability.
This is a calculation of sheaves over $F$: in local symplectic
eigenframes the completed twisted Koszul contractions are
$\mathcal O_F$-linear. The moving directions use $(1-g)^{-1}$
and $\hbar$-adic perturbation, and the fixed directions use
the formal Poincar\'e contraction divided by $\hbar$.
Thus the calculation remains valid over the base ring, and
\et\ matrix splittings give the line $\mathscr L_g((\hbar))$
in degree $d_F$.
To identify the induced connection, locally write
$g=\operatorname{Ad}(u)$ on $\mathscr D|_F$ and take $u^q=1$
after a finite \et\ extension. The twisted matrix trace is
$ag\mapsto\operatorname{tr}(au)$, with generator $u^{-1}g$.
Changing $u$ by a $q$-th root of unity gives exactly the
transition functions of $\mathscr L_g$. Invariance of the
algebra connection gives $u^{-1}\nabla^{\mathscr D}u=\alpha\,1$;
differentiating $u^q=1$ gives $q\alpha=0$. Thus this generator
is horizontal for the canonical finite connection. The Weyl
generator, normalized by the symplectic volume of $T_F$, is
also horizontal for the product connection $\nabla$.

Filter by exterior degree before taking derived global sections.
The fiber Hochschild homology sheaf of its $p$-th graded piece
is $\Omega_F^p\otimes\mathscr L_g((\hbar))$ in chain degree
$d_F$. Put $\mathcal E=\Omega_F^p\otimes\mathscr L_g$.
Its finite parameter quotients are quasi-coherent with
surjective transition maps, so
$\mathcal E[[\hbar]]=R\!\varprojlim_n\mathcal E[\hbar]/\hbar^n$
has no higher cohomology on affine opens
\cite[Tag~0A0J]{StacksReflexive}.
Localization in $\hbar$ is a filtered colimit, which commutes
with sheaf cohomology on the Noetherian $F$
\cite[Tag~01FF]{StacksReflexive}. Moreover, $F$ has finite
Zariski cohomological dimension
\cite[Tag~02UZ]{StacksReflexive}; hence localization also
commutes with derived global sections of the bar total complex,
despite its unbounded chain degrees. The first page is
\[
 E_1^{p,-j}=
 \begin{cases}
 \Gamma(F,\Omega_F^p\otimes\mathscr L_g)((\hbar)),&j=d_F,\\
 0,&j\ne d_F.
 \end{cases}
\]
The connection $\nabla$ commutes with the fiber Hochschild
differential. The remaining terms of $D$ are inner derivations
with $g$-invariant Hamiltonians along $F$, and their action on
twisted Hochschild homology is zero by the insertion homotopy.
For $\delta|_F$, invariance follows from the symplectic
orthogonality of the fixed and moving tangent spaces.
Consequently $d_1=\nabla_{\mathscr L_g}$.
Since $0\le p\le d_F$, the spectral sequence converges;
its single nonzero row excludes all higher differentials and
leaves a single filtration quotient in each total degree.
It therefore identifies the total cohomology with
\[
 H^{d_F-i}_{\mathrm{dR}}(F,\mathscr L_g)((\hbar))
\]
in Hochschild degree $i$. For scalar coefficients this is the Fedosov comparison
of~\cite[Proposition~12 and (69)]{DE}.
The connection on $\mathscr L_g$ is trivialized by a finite \et\ cover, hence is regular
singular. Deligne's comparison~\cite[II, Th\'eor\`eme~6.2]{DeligneRegular}
gives~\eqref{qr:twisted-hh}. Its finite-dimensional classical
cohomology commutes with the formal coefficient extension.
Every comparison is equivariant and commutes with restriction
to stable opens.

For the root realization, matrix duality identifies the
$\chi$-summand of the original family with the $g$-twisted
summand, where $g$ is paired with $\chi$ by
\eqref{qm:sector-label}. Consequently
Theorem~\ref{ql:completion} gives
\begin{equation}\label{qr:sector-realization}
 \HH_i(\mathcal A)_\chi\otimes_{E(t)}\C((\hbar))
 \simeq
 \bigoplus_a H^{N-2a-i}
 \bigl(((X^\vee)^g)_a,\mathscr L_{g,\hbar}\bigr)^G
\end{equation}
at $t=\rho e^\hbar$. At $t=e^\hbar$ the same generic summand
is $H^{N-i}(X)_\chi((\hbar))$.
These statements also hold on every simultaneous source and
dual freezing chart. Indeed, the denominators are homogeneous,
inverting a central Frobenius coordinate in the crossed product
is equivalent to inverting the corresponding source variable,
and the crossed products,
matrix equivalence, and natural bar-chain comparisons all commute
with these localizations. In applying this assertion we retain
the inherited groups $Q,G$ and the inherited diagonal matrix $T$,
including its entries at newly frozen variables.

\begin{lemma}\label{qr:isolated-lowest}
Let $X_\nu$ and $Y_\nu\subset X^\vee$ be corresponding isolated
freezing charts, and let $g\in G$ have nonempty fixed locus in
$Y_\nu$. Write $S$ for its moving mutable set, $a=|S|$, and
$\chi$ for its paired source character. Then
\[
 \dim H^{2a}(X_\nu)_\chi=1,
 \qquad \dim H^0(Y_\nu^g)^G=1.
\]
Every component of $Y_\nu^g$ has age $a$.
\end{lemma}

\begin{proof}
Write $\alpha_i$ for the retained source exchange columns and
$n_\nu$ for their number. The component characters of the
diagonal automorphism group are
\[
 \Xi_\nu=(\operatorname{span}_{\Q}\{\alpha_i\}\cap\Z^N)
               /\sum_i\Z\alpha_i.
\]
The finite pairing of Lemma~\ref{qm:finite-pairing} remains
injective on this group after freezing: it sends
$[\sum_i r_i\alpha_i]$ to
$w\mapsto\exp(-2\pi i\sum_i r_iw_i)$, whose kernel consists
exactly of integral $r_i$. Under this identification, the support
$\{i:r_i\notin\Z\}$ of $\chi$ is the moving set $S$ of $g$.
Lam--Speyer's isolated character formula
\cite[Theorem~7.5]{LS} therefore gives
\[
 \sum_k\dim H^k(X_\nu)_\chi\,z^k
 =z^{2a}(1+z+z^2)^{n_\nu-a}(1+z)^{N-2n_\nu}.
\]
Its lowest coefficient proves the first equality. The theorem
applies directly to the zero principal block of an isolated
exchange matrix.

On the dual chart, the equations $x_i=y_i=0$ for $i\in S$
cut out a frozen-torus fiber, together with the remaining isolated
factors. Its components are smooth and irreducible: on each
connected torus fiber the remaining exchange characters are
independent, and the usual isolated presentation embeds in its
Laurent torus. The component group is dual to the saturation
quotient of the columns indexed by $S$.

To determine the action of $G$ on these components, write $J_i$
for the columns of the current-seed matrix $J$. The retained dual
columns are its mutable columns, and the character group of $G$
is $\Z^N/J\Z^N$. There are injections
\[
 \frac{\operatorname{span}_{\Q}\{J_i:i\in S\}\cap\Z^N}
      {\sum_{i\in S}\Z J_i}
 \lhook\joinrel\longrightarrow
 \frac{\operatorname{span}_{\Q}\{J_i:i\text{ mutable}\}\cap\Z^N}
      {\sum_{i\text{ mutable}}\Z J_i}
 \lhook\joinrel\longrightarrow
 \frac{\Z^N}{J\Z^N}.
\]
Indeed, applying $J^{-1}$ reduces both assertions to the
corresponding coordinate sublattices. Dualizing shows that $G$
surjects onto the component group of the diagonal grading torus
after freezing and acts transitively on the components of $Y_\nu^g$.
Hence $\dim H^0(Y_\nu^g)^G=1$.

Each moving pair has eigenvalues $h_i,h_i^{-1}$ and contributes
one to the age, proving the last assertion.
\end{proof}

\begin{lemma}\label{qr:untwisting}
For every $g\in G$ with nonempty fixed locus, the trace local
system $\mathscr L_{g,\hbar}$ on $(X^\vee)^g$ is
$G$-equivariantly trivial on each orbit of connected components.
\end{lemma}

\begin{proof}
Choose an isolated freezing chart $Y_\nu$ meeting a fixed
component, and use the notation of Lemma~\ref{qr:isolated-lowest}.
The two realizations~\eqref{qr:sector-realization} of the same
generic character summand, in Hochschild degree $N-2a$, give
\[
 \dim_{\C((\hbar))}
 H^0(Y_\nu^g,\mathscr L_{g,\hbar})^G
 =\dim_{\C}H^{2a}(X_\nu)_\chi=1.
\]
A rank-one local system contributes at most one invariant
section on an orbit of components. There is exactly one such
orbit, and equality forces both its monodromy and the action
of a component stabilizer on its flat sections to be trivial.

A connected component $F$ of $(X^\vee)^g$ is smooth and
irreducible. Its intersection with $Y_\nu$ is a nonempty
Zariski open set, whose fundamental group surjects onto
$\pi_1(F)$. Thus the flat trivialization extends to $F$.
The chart is $G$-stable, so the stabilizer of $F$ preserves this
open set and acts trivially on the extended section as well.
Transporting it among the components proves the assertion.
\end{proof}

\subsection{Periods and the cap action}

Write the period of the model in Lemma~\ref{qr:fedosov-model} as
\begin{equation}\label{qr:period}
 \operatorname{per}(\mathscr D_\hbar)
 =\hbar^{-1}[\Omega_\hbar]
 =\hbar^{-1}[\omega_Y]+\eta_0+\hbar\eta_1+\hbar^2\eta_2+\cdots,
 \qquad \eta_r\in H^2_{\mathrm{dR}}(Y)^\Pi.
\end{equation}

\begin{proposition}\label{qr:cap-period}
Under \eqref{qr:twisted-hh}, cap product with the actual
Kodaira--Spencer class acts on the summand indexed by $F$ as cup
product with
\[
 c\,\left.\partial_\hbar
          \operatorname{per}(\mathscr D_\hbar)\right|_F,
 \qquad c\in\C^\times.
\]
The constant depends only on the normalization of the Hochschild
and star-product conventions and is the same for every sector.
The formula is equivariant and remains valid after transporting the
comparison along a parameter-dependent equivalence of formal deformations.
\end{proposition}

\begin{proof}
We first specify the cochain and chain comparison. Let
$\mathcal H^{p,q}=\Omega_Y^p(C^q(\mathcal W_{\mathscr D}))$,
with total degree $p+q$ and differential
$d_{\mathrm{tot}}=D+(-1)^p b$.
Here $C^q$ consists of normalized
$\mathcal O_Y((\hbar))$-multilinear fiber cochains, continuous
for the combined Fedosov filtration and with a uniform lower
bound on their $\hbar$-orders, as in
\cite[Definition~1 and Appendix~B]{Dolgushev}.
We take derived
global sections throughout. Central zero-cochains give a morphism
of differential graded algebras
\[
 j:R\Gamma(Y,\Omega_Y^\bullet((\hbar)))\longrightarrow
                 R\Gamma(Y,\operatorname{Tot}\mathcal H^{\bullet,\bullet}).
\]
The bounded exterior-degree filtration proves this is a
quasi-isomorphism: fiberwise Weyl cohomology is the center in
degree zero, and the Azumaya factor is separable.
The flat-lift comparison with the differential cochains of
$\mathscr D_\hbar$ is the Azumaya-valued version of
\cite[Theorem~2 and Proposition~6]{Dolgushev}.
For the flat lifts $\tau,\tau_g$ of algebra and module sections,
and the lift $P$ of a cochain $P_0$, the evaluation identities are
\cite[(3.49)]{Dolgushev}
\[
 \tau(P_0(a_1,\ldots,a_q))
   =P(\tau(a_1),\ldots,\tau(a_q)),\qquad
 \tau_g(mg(a))=\tau_g(m)g(\tau(a)).
\]
Substitution in the bar formula proves compatibility with cap
product. On fiber homology the central action also commutes
with matrix trace.
For a closed form $\alpha$, the action of $j(\alpha)$ on the
exterior-degree spectral sequence is multiplication by
$\alpha|_F$. On its first page this is the wedge action on
$\Omega_F^\bullet\otimes\mathscr L_g$. The single-row
identification above therefore carries its action on total
cohomology to cup product by $[\alpha|_F]$.

We calculate the class globally. Write
$D=\nabla+\operatorname{ad}(\gamma/\hbar)$, incorporating
$-\delta$ into the linear part of $\gamma$. Its central curvature
equation is
\begin{equation}\label{qr:central-curvature}
 C=R_W/\hbar+F_{\mathscr D}
          +\nabla(\gamma/\hbar)+(\gamma/\hbar)^2
       =-\Omega_\hbar/\hbar.
\end{equation}
Let $\mathsf E$ be the Euler operator in the Weyl fiber coordinates.
The operator
\[
 T=\partial_\hbar+\frac{1}{2\hbar}\mathsf E
\]
is a derivation of the Weyl product, extending $\partial_\hbar$
on scalars, and acts trivially on the classical Azumaya factor.
Indeed, a term of order $\hbar^j$ in that product lowers fiber
degree by $2j$, so the two multiplication defects cancel.
The connections are independent of $\hbar$, and hence
$[T,\nabla]=0$. Moreover, $R_W$ is quadratic in fiber variables,
so $T(R_W/\hbar+F_{\mathscr D})=0$. Put
\[
 v=T(\gamma/\hbar)\in
             \Omega_Y^1(\mathcal W_{\mathscr D})((\hbar)).
\]
Differentiating~\eqref{qr:central-curvature} gives the two global
identities
\begin{equation}\label{qr:global-ks-curvature}
 [T,D]=\operatorname{ad}(v),\qquad
 Dv=-\partial_\hbar(\Omega_\hbar/\hbar).
\end{equation}

We identify this class with the actual parameter class using
the same flat-lift comparison. For the chosen Fedosov product,
prolonging the flat lift to jets
gives a global filtered algebra isomorphism
\[
 \lambda_\hbar:\mathscr J(\mathscr D)[[\hbar]]
     \longrightarrow\mathcal W_{\mathscr D}
\]
where $\mathscr J$ denotes formal jets along the diagonal.
The source has its $\mathcal O_Y$-module structure from the
first projection and the prolonged star product in the second
factor. The map is $\mathcal O_Y[[\hbar]]$-linear and intertwines
the Grothendieck connection and $D$.
Modulo $\hbar$ it is the filtered Taylor isomorphism; at each
parameter order and finite jet order the construction uses
only finitely many jets. The $\hbar$-adic lifting therefore
gives this continuous isomorphism and its inverse. Transport
the coefficientwise parameter connection
on the jet bundle to $T_0=\lambda_\hbar\partial_\hbar^{\mathscr J}
\lambda_\hbar^{-1}$. Then $[T_0,D]=0$, and its multiplication
defect is the flat lift $\widetilde\kappa$ of the actual
Kodaira--Spencer cocycle. Both operators extend
$\partial_\hbar$ on scalars and kill the unit, so
$P=T-T_0$ is $\mathcal O_Y((\hbar))$-linear and normalized.
The maps $\lambda_\hbar^{\pm1}$ preserve the integral
$\hbar$-lattices and are continuous for the combined filtration;
thus $P$ is continuous and lowers $\hbar$-valuation by at most
one. In particular, it belongs to the specified fiber cochain
complex. Since $T$ is a derivation,
\[
 bP=\widetilde\kappa,\qquad
 DP=-[T,D],\qquad
 d_{\mathrm{tot}}P=\widetilde\kappa-[T,D].
\]
Thus the actual class is represented by $[T,D]$ in the
fiberwise double complex.

For the zero-cochain
$v\in\mathcal H^{1,0}$, our convention $bv(a)=av-va$ gives
\[
 d_{\mathrm{tot}}v=Dv+\operatorname{ad}(v).
\]
In particular,
\[
 \widetilde\kappa
   -j\bigl(\partial_\hbar(\Omega_\hbar/\hbar)\bigr)
       =d_{\mathrm{tot}}(P+v).
\]
The cap differential identity gives a global homotopy of total
chain maps. This homotopy need not preserve the exterior
filtration: we first replace the cap action by the central
representative, and only then use its filtration-degree-two
action on the spectral sequence above.
Thus, allowing for the normalization of the Weyl product,
its image under $j^{-1}$ is
$c\,\partial_\hbar[\Omega_\hbar/\hbar]$.
This is the characteristic-class calculation of
\cite[Theorem~4.1]{WX}, now performed in the single global
Azumaya-valued resolution. The constant $c\ne0$ accounts only
for the Hochschild and Weyl conventions. In the conventions
of that reference,
$d_{\mathrm{WX}}=(i/\hbar)d_{\mathrm{Hoch}}$ and cochains of
degree $q$ are multiplied by $(i/\hbar)^q$; this converts
their factor $-i\hbar^2$ in degree two into a constant.
No sector has entered this computation.

The cap comparison gives the asserted action. Invariance of the
chosen connections and curvatures permits finite invariants;
Lemma~\ref{qr:ks-completion} preserves the class upon returning
to the original formal product.
\end{proof}

Apply the proposition first to $(Y,\Pi)=(X,Q)$ and the ordinary
star product, and then to $(Y,\Pi)=(X^\vee,G)$ and $\mathscr B_\hbar$.
By \eqref{qr:second-purity}, all coefficients of the corresponding
periods have type $(2,2)$. Multiplying each cap operator by a
nonzero Laurent scalar gives cup product with, respectively,
\begin{equation}\label{qr:deformed-lefschetz}
 \ell_\hbar=[\omega]+O(\hbar)
       \in H^{2,(2,2)}(X)^Q[[\hbar]],\qquad
 \ell_\hbar^\vee=[\omega^\vee]+O(\hbar)
       \in H^{2,(2,2)}(X^\vee)^G[[\hbar]].
\end{equation}
The second acts by restriction on every dual fixed sector.
Lemma~\ref{qr:untwisting} trivializes its trace local system, and
matrix duality identifies the finite source character.

\subsection{Recovering the age-shifted Tate degrees}

\begin{lemma}\label{qr:monodromy}
Let $Z$ be smooth of dimension $2e$, with split mixed Tate
cohomology, and let $\omega_Z\in H^{2,(2,2)}(Z)$ satisfy
curious Lefschetz.  Suppose
$\ell=[\omega_Z]+O(\hbar)$ has all its coefficients in
$H^{2,(2,2)}(Z)$.  On $H^*(Z)((\hbar))$, the monodromy
filtration centered at zero of cup product with $\ell$ is
\begin{equation}\label{qr:monodromy-formula}
 M_j=\bigoplus_{p\ge e-j} H^{*,(p,p)}(Z)((\hbar)).
\end{equation}
The same assertion holds on an eigenspace for a finite group
preserving $\ell$ and $\omega_Z$.
\end{lemma}

\begin{proof}
By Deligne's characterization
\cite[Proposition~1.6.1]{DeligneWeilII}, it suffices to check
$\ell M_j\subseteq M_{j-2}$ and that
$\ell^j:\operatorname{Gr}_j^M\to\operatorname{Gr}_{-j}^M$
is an isomorphism for $j\ge0$. The first follows from the Tate
degree of $\ell$; the second follows because its determinant
has the nonzero constant term supplied by curious Lefschetz.
The same argument applies to each finite-group eigenspace.
\end{proof}

Let $\chi\in S_A^\vee$, regarded as a character of $Q$, with
paired element $g$ given by \eqref{qm:sector-label}, and put
\[
 V_i=\HH_i(\mathcal A)_\chi,
 \qquad L=\kappa\cap-.
\]
Theorem~\ref{ql:completion} and the realization at
$t=e^{\hbar}$ show that $V=\bigoplus_iV_i$ is finite
dimensional and that $L$ is nilpotent.  Let $M$ denote its
monodromy filtration centered at zero.  It is homogeneous for
the Hochschild grading and commutes with field extension,
since it is determined by the kernels and images of powers
of $L$.  Multiplication of $L$ by a nonzero scalar preserves $M$.

\begin{theorem}\label{qr:character-comparison}
For every $k,p$ and every $\chi\in S_A^\vee$ with paired element $g$,
\begin{equation}\label{qr:hodge-identity}
 \dim H^{k,(p,p)}(X)_\chi
 =\sum_a
   \dim H^{k-2a,(p-a,p-a)}\bigl(((X^\vee)^g)_a\bigr)^G,
\end{equation}
where $((X^\vee)^g)_a$ is the union of dual fixed components of
age $a$.
\end{theorem}

\begin{proof}
At $t=e^{\hbar}$, the continuous Hochschild realization is
$H^{N-i}(X)_\chi((\hbar))$.  Proposition~\ref{qr:cap-period}
and \eqref{qr:deformed-lefschetz} identify $L$, up to a
nonzero scalar, with cup product by $\ell_\hbar$.
Proposition~\ref{geom:fixed-lefschetz} with $h=1$ and
Lemma~\ref{qr:monodromy} therefore give
\begin{equation}\label{qr:source-multiplicity}
 \dim H^{k,(p,p)}(X)_\chi
 =\dim_{E(t)}
       \operatorname{Gr}_{d-p}^M V_{N-k}.
\end{equation}

At $t=\rho e^{\hbar}$, Theorem~\ref{qm:morita}
and Theorem~\ref{ql:completion} instead identify the same
graded vector space with
\[
 \bigoplus_a
 H^{N-2a-i}\bigl(((X^\vee)^g)_a\bigr)^G((\hbar)).
\]
Indeed, Lemma~\ref{qr:untwisting} gives a $G$-equivariant flat
trivialization of the trace local system, and
\eqref{qr:sector-realization} identifies its sector with precisely
the character $\chi$. The trivialization commutes with cup product by every restricted period
coefficient. Proposition~\ref{qr:cap-period} therefore identifies
$L$ on this realization, up to a nonzero Laurent scalar, with cup
product by $\ell_\hbar^\vee$.

A component $F$ of age $a$ is symplectic of dimension
$N-2a=2(d-a)$.  By Proposition~\ref{geom:fixed-lefschetz},
$\omega^\vee|_F$ satisfies curious Lefschetz. Its monodromy
filtration is thus given by
\[
 p_F\ge d-a-j
 \quad\Longleftrightarrow\quad
 p_F+a\ge d-j.
\]
In Hochschild degree $i=N-k$, its cohomological degree is
$k-2a$.  Hence $\operatorname{Gr}_{d-p}^M V_{N-k}$ in
this second realization has dimension
\[
 \sum_a\dim
 H^{k-2a,(p-a,p-a)}\bigl(((X^\vee)^g)_a\bigr)^G.
\]
Comparison with \eqref{qr:source-multiplicity} proves
\eqref{qr:hodge-identity}.
\end{proof}

\section{Motivic mirror symmetry}
\label{sec:main-proof}

Each component of $T_A$ contains a torsion point. Two torsion lifts of the same
element of $S_A$ differ by a torsion point of $T_A^0$, which acts
trivially on $M(X)$ by Lemma~\ref{geom:torsion}. Thus torsion lifts
define an action of $S_A$ on $M(X)$ independent of their choice.
For a character $\chi$ of $S_A$ and a number field $E$ splitting
$S_A$, put
\[
 e_\chi=\frac{1}{|S_A|}\sum_{s\in S_A}\chi(s)^{-1}M(\widetilde s),
 \qquad M(X;E)_\chi=e_\chi M(X;E),
\]
where $\widetilde s$ is any torsion lift. On cohomology this is the
projector for the pullback action used in
Theorem~\ref{qr:character-comparison}.

\begin{theorem}\label{main:character}
Let $A$ be full-rank and sink-recurrent, with skew-symmetrizable principal part, and let
$G=G_0$ be the group of Proposition~\ref{geom:mirror}.
If $\chi=[Ar]\in\sat(A\Z^n)/A\Z^n$ and
$h=(\exp(2\pi i r),1)$, then
\begin{equation}\label{main:character-motive}
 M(X;E)_\chi\simeq
 \bigoplus_a M((X^\vee)^h_a;E)^G(a)[2a].
\end{equation}
In even total dimension, one may also take $G=G_J$, where
$J=P^{-1}K$ for any nonsingular integral alternating completion $K$ of
$A\operatorname{diag}(p_i)$.
In every total dimension there are rational isomorphisms
\begin{equation}\label{main:motives}
 M(X)\simeq\Morb(\cM_A),\qquad
 M^c(X)\simeq M^c_{\mathrm{orb}}(\cM_A),
\end{equation}
where the compact orbifold motive uses the same age twists.
\end{theorem}
\begin{proof}
Suppose first that $N$ is even and $G=G_J$.
By Theorem~\ref{geom:splitting}, both sides of
\eqref{main:character-motive} are finite split Tate motives with
$E$ coefficients: for the character summand, apply
Lemma~\ref{geom:arithmetic-splitting} to the corresponding mixed Tate
retract over a sufficiently large number field and its split Hodge realizations.
The multiplicity of $E(p)[k]$ on its right side is
\[
 \sum_a h^{k-2a,p-a}((X^\vee)^h_a)^G,
\]
which is the source multiplicity by~\eqref{qr:hodge-identity}.
Choosing bases of the multiplicity spaces gives the noncanonical
isomorphism~\eqref{main:character-motive}.

The source projectors sum to the identity. Conversely all elements of $G$ with
nonempty fixed loci are the elements~\eqref{qm:sector-label}.
Proposition~\ref{geom:inertia-comparison} replaces $G_J$ by $G_0$
in each fixed motive, proving the character formula for $G_0$ as well.
Summing the character multiplicities proves the rational identity: both
rational motives are split Tate, and equality after the faithful
coefficient extension $\Q\subset E$ identifies their rational
multiplicities.

For odd $N$, adjoining an unused frozen coordinate multiplies the source,
mirror, and fixed loci by $\Gm$, preserving $S_A$, $G_0$, their pairing,
and the ages. The even-dimensional comparison for $G_0$ equates their characterwise Tate
multiplicity polynomials after multiplication by $1+uv$, where $u$
records cohomological degree and $v$ Tate degree. Cancellation in
$\Z[u,v]$ and Theorem~\ref{geom:splitting} give the result for odd $N$.

Finally, smooth duality in dimension $N$ transforms a fixed term into
\[
 \bigl(M(F)(a)[2a]\bigr)^\vee(N)[2N]
     =M^c(F)(a)[2a],\qquad \dim F=N-2a.
\]
It commutes with finite invariants, proving the compact assertion.
\end{proof}

\begin{corollary}\label{main:chow-k}
Let $Z$ be a smooth quasi-projective complex variety. There are
additive rational isomorphisms
\begin{align*}
 \CH^p(X\times Z,r)_\Q
 &\simeq\bigoplus_{g,a}
       \CH^{p-a}((X^\vee)^g_a\times Z,r)_\Q^G,\\
 K_j(X\times Z)_\Q
 &\simeq\bigoplus_{g,a}K_j((X^\vee)^g_a\times Z)_\Q^G,
 \qquad r,j\ge0,
\end{align*}
where $G$ is the finite group of Proposition~\ref{geom:mirror}.
In the $K$-theory formula, Adams degree $p$ on the source corresponds
to Adams degree $p-a$ on a sector of age $a$.
\end{corollary}
\begin{proof}
Tensor~\eqref{main:motives} with $M(Z)$ and apply
$\Hom(-,\Q(p)[2p-r])$ to obtain the higher Chow comparison.
The rational cycle-class comparison identifies
$\CH^p(-,j)_\Q$ with the Adams summand $K_j^{(p)}$ by
\cite[Theorems~3.1 and~4.7]{LevineChow}. This is also the rational
Adams decomposition of the motivic spectral sequence of~\cite{FS}.
Summing over $p$ proves the second assertion.
\end{proof}

For a finite \et\ $G$-torsor $p:V\to W$ of degree $d$ between
smooth finite-type varieties, the graph and transpose-graph
correspondences give
\begin{equation}\label{main:transfer-identities}
 M(p)\operatorname{Tr}_p=d\,\id,\qquad
 \operatorname{Tr}_pM(p)=\sum_{g\in G}M(g).
\end{equation}
The second identity uses $V\times_WV=\coprod_g\Gamma_g$.
When $B$ is skew-symmetric and $A$ is really full-rank, $T_A$ is connected, so
Lemma~\ref{geom:torsion} makes every $g\in G_0$ act trivially on
$M(X)$ after inverting $d$. Applying~\eqref{main:transfer-identities}
to the morphism $p_A$ of Proposition~\ref{geom:rfr} gives
\[
 M(p_A)^{-1}=d^{-1}\operatorname{Tr}_{p_A}
 \qquad\text{with coefficients }\Z[1/d].
\]
Thus in this case the motivic comparison is realized by the specified
finite \et\ morphism, and its pullback on higher Chow groups and
algebraic $K$-theory also preserves products and natural operations.

\section{The projective-line interpretation}\label{sec:p1}

Assume throughout this section that $B$ is skew-symmetric and that
$A$ is really full-rank and sink-recurrent.
Set $Y=\cM_A$ and let $p=p_A:X\to Y$ be the
finite \et\ $G$-torsor of Proposition~\ref{geom:rfr}, of degree $d$.
Deformation to the normal cone extends this comparison to
torus-stable smooth closed pairs over $\mathbb P^1$.

Put $S=\mathbb P^1_\C$, write $\pi:S\to\Spec\C$, and use the
coordinate $t$ on $S\setminus\{\infty\}$ and parameter $t^{-1}$ at $\infty$.
In the stable $\infty$-category $\mathrm{DM}(S;\Q)$, define
$M_S(V)=f_\sharp\Q_V$ for a smooth $f:V\to S$, and
\[
 M_S(V,V^\circ)=\operatorname{cofib}
   \bigl(M_S(V^\circ)\longrightarrow M_S(V)\bigr)
\]
for an open $V^\circ\subset V$. These constructions commute with
base change.

\subsection{Deformation of a smooth pair}

Let $Z\hookrightarrow V$ be a closed immersion of smooth complex
varieties of pure codimension $a\geq1$, and put $V^\circ=V\setminus Z$.
The deformation to the normal cone over $S$ is
\[
 \mathcal D_ZV=
 \operatorname{Bl}_{Z\times\{\infty\}}(V\times S)
       \setminus\widetilde{V\times\{\infty\}},
 \qquad Z_S=Z\times S\subseteq\mathcal D_ZV,
\]
where the tilde denotes the strict transform. It is smooth over $S$,
and the fibres of
$(\mathcal D_ZV,\mathcal D_ZV\setminus Z_S)$ at $0$ and $\infty$
are respectively
\[
 (V,V^\circ),\qquad (N_{Z/V},N_{Z/V}\setminus Z),
\]
as in~\cite[Section~2.3.1]{Deglise2018}. Homotopy purity
\cite[Section~3.2, Theorem~2.23]{MV99}, in the relative form of
\cite[Theorem~2.3.2]{Deglise2018}, and the oriented Thom
isomorphism \cite[(2.2.22.a)]{Deglise2018} give
\begin{equation}\label{eq:relative-purity-p1}
 \mathcal C_{Z/V}:=
 M_S(\mathcal D_ZV,\mathcal D_ZV\setminus Z_S)
       \simeq\pi^*M(Z)(a)[2a].
\end{equation}
Indeed, the normal bundle of $Z_S$ is
$\operatorname{pr}_Z^*N_{Z/V}\otimes\operatorname{pr}_S^*\mathcal O_S(-\infty)$,
of rank $a$, and $M_S(Z_S)=\pi^*M(Z)$.
By naturality, the two fibre pullbacks of~\eqref{eq:relative-purity-p1}
are the canonical purity identifications for $(V,Z)$ and
$(N_{Z/V},Z)$. The equivalence also respects finite-group actions
on the closed pair that fix $S$. Composing the identification at $0$
with $M(V)\to M(V,V^\circ)$ gives the Gysin map for
$Z\hookrightarrow V$. Passing to the ordinary motive of the total
pair gives
\[
 M(\mathcal D_ZV,\mathcal D_ZV\setminus Z_S)
       \simeq M(Z)(a)[2a]\otimes M(\mathbb P^1).
\]

\subsection{Compatibility with the mirror quotient}

\begin{proposition}\label{prop:mirror-dnc}
Let $Z\subseteq X$ be a nonempty smooth closed $T_A$-stable
subvariety of pure codimension $a\geq1$, and set $W=p(Z)\subseteq Y$.
Then $W$ is smooth, $Z=p^{-1}(W)$, and $p_Z:Z\to W$ is a finite
\et\ $G$-torsor of degree $d$.  The deformation spaces fit into a
Cartesian square
\[
 \begin{tikzcd}
 \mathcal D_ZX \arrow[r,"p_{\mathcal D}"] \arrow[d]
   & \mathcal D_WY \arrow[d] \\
 X \arrow[r,"p_A"] & Y,
 \end{tikzcd}
\]
and $p_{\mathcal D}$ is a finite \et\ $G$-torsor over $S$.
At $0$ it restricts to $p_A$; at $\infty$ it restricts to the
normal-bundle morphism induced by
$N_{Z/X}\simeq p_Z^*N_{W/Y}$.

The induced morphism of relative motives is an isomorphism, and
under purity it is the bottom arrow in the commutative diagram
\begin{equation}\label{eq:mirror-dnc-purity}
 \begin{tikzcd}[column sep=7em]
 \mathcal C_{Z/X}
    \arrow[r,"M_S(p_{\mathcal D})"] \arrow[d,"\sim"']
   & \mathcal C_{W/Y} \arrow[d,"\sim"] \\
 \pi^*M(Z)(a)[2a]
    \arrow[r,"{\pi^*M(p_Z)(a)[2a]}"']
   & \pi^*M(W)(a)[2a].
 \end{tikzcd}
\end{equation}
Its inverse, under these identifications, is
\[
 \pi^*\!\left(\frac1d M({}^{\mathsf T}\Gamma_{p_Z})\right)(a)[2a].
\]
\end{proposition}

\begin{proof}
The $G$-stable closed subvariety $Z$ descends along the torsor $p$ to
a closed subvariety $W\subseteq Y$, with $Z=X\times_YW$.
Smoothness and codimension descend along $p_Z$. Flat base change
for blowups and removal of strict transforms give the Cartesian
square and the normal-bundle identification. Thus $p_{\mathcal D}$
is a finite \et\ torsor of degree $d$.

Purity is compatible with this \et\ Cartesian square.  Together
with the naturality of the Thom isomorphism, this gives
\eqref{eq:mirror-dnc-purity}.  Since $Z$ is $T_A$-stable, the action
of the connected torus $T_A$ restricts to $Z$.  Lemma~\ref{geom:torsion}
therefore gives $M(g)=\id$ on $M(Z)$ for every $g\in G$.
The transfer identities~\eqref{main:transfer-identities}, applied
to $p_Z$, give the stated inverse.
Pullback to $0$ and $\infty$ preserves the purity diagram and its
inverse, proving the endpoint assertions.
\end{proof}

In particular, the comparison respects the Gysin maps:
\[
 \begin{tikzcd}[column sep=7em]
 M(X) \arrow[r,"M(p_A)"] \arrow[d,"\operatorname{Gy}_{Z/X}"']
   & M(Y) \arrow[d,"\operatorname{Gy}_{W/Y}"] \\
 M(Z)(a)[2a] \arrow[r,"{M(p_Z)(a)[2a]}"']
   & M(W)(a)[2a].
 \end{tikzcd}
\]

\appendix
\section{Hochschild homology of isolated quantum cluster algebras}
\label{qt:section}

We prove the isolated-seed comparison with continuous Hochschild
homology. Let $E$ be a field of characteristic zero; all commutation parameters
are powers of $t=q^{1/2}$.

Let
\[
 L=L_f\oplus\bigoplus_{i=1}^n\Z e_i
\]
be a lattice, and let $\Lambda$ be an integral alternating form on $L$
which is nondegenerate on $L_{\Q}$. Suppose that
$c_1,\ldots,c_n\in L_f$ are linearly
independent and satisfy
\begin{equation}\label{qt:compatibility}
 \Lambda(c_i,L_f)=0,
 \qquad
 \Lambda(c_i,e_j)=d_i\delta_{ij},
 \qquad d_i\ne0.
\end{equation}
Put $R_0=E[t^{\pm1}]$ and
$F=E(t)$. The quantum torus associated with $L$ has basis
$\{Z^v:v\in L\}$ and multiplication
\[
 Z^v Z^w=t^{\Lambda(v,w)}Z^{v+w}.
\]
Choose $\eta_i\in E^\times t^{\Z}$, and let $\mathcal A$ be its
$R_0$-subalgebra generated by the frozen monomials $Z^f$, $f\in L_f$,
and
\begin{equation}\label{qt:generators}
 x_i=Z^{e_i},
 \qquad
 y_i=(1+\eta_i Z^{c_i})x_i^{-1}
 \qquad(1\le i\le n).
\end{equation}
This normalizes an isolated quantum seed by multiplying its exchange
partners by frozen units. Write
$A=\mathcal A\otimes_{R_0}F$.

Fix $t_0\in E^\times$ and use the embedding
\[
 R_0\longrightarrow E[[\hbar]],
 \qquad t\longmapsto t_0\exp(\hbar).
\]
Set
\[
 K=E((\hbar)),\qquad
 \widehat A_0=\varprojlim_{r}
 \bigl(\mathcal A\otimes_{R_0}E[[\hbar]]\bigr)/\hbar^r,
 \qquad
 \widehat A=\widehat A_0[\hbar^{-1}].
\]
The continuous Hochschild complex of \eqref{ql:continuous-definition},
after inverting $\hbar$, has chain modules
\begin{equation}\label{qt:continuous-chains}
 C_j^{\mathrm{cont}}(\widehat A/K)
 =\left(\varprojlim_r
   (\widehat A_0/\hbar^r)^{\otimes_{E[[\hbar]]/\hbar^r}(j+1)}
   \right)[\hbar^{-1}],
\end{equation}
with the Hochschild differential. Their homology is denoted by
$\HH_j^{\mathrm{cont}}(\widehat A/K)$.

\begin{theorem}\label{qt:comparison}\label{qt:completion}
For the isolated quantum algebra defined by
\eqref{qt:compatibility} and \eqref{qt:generators}, the canonical map
\[
 C_\bullet(A/F)\otimes_F K
 \longrightarrow C_\bullet^{\mathrm{cont}}(\widehat A/K)
\]
is a quasi-isomorphism.
\end{theorem}

\subsection{Completion and inner conjugation}

The standard isolated-seed monomials form an $R_0$-basis of
$\mathcal A$ \cite[proof of Proposition~8.2]{MullerQuantum}:
in a fixed ordering, they consist of a frozen monomial and, for each
$i$, a nonnegative power of $x_i$ or a positive power of $y_i$,
but not both.

For any such basis $\{b_\alpha\}$, the completed generic fiber consists
of sums
\begin{equation}\label{qt:restricted-support}
 \sum_\alpha a_\alpha b_\alpha,
 \qquad a_\alpha\in K,
\end{equation}
whose $\hbar$-adic valuations are bounded below and for which, for each
integer $M$, only finitely many coefficients have valuation less than
$M$. This description applies to the tensor bases in
\eqref{qt:continuous-chains} as well.

A linear map with finite images of basis vectors and a uniform lower
valuation bound $-c$ extends to~\eqref{qt:restricted-support}: its output
modulo $\hbar^M$ depends only on the finitely supported input modulo
$\hbar^{M+c}$. No bound on the lengths or exponents of those images
is needed. We apply this criterion degreewise to chain maps and homotopies.

Let $B$ be an algebra, let $g$ be an automorphism, and write $B_g$ for
the bimodule with its usual left action and right action
$b\cdot a=b\,g(a)$. An invertible element $u\in B$ acts on
$C_\bullet(B,B_g)$ by
\begin{equation}\label{qt:inner-action}
 T_u(m\mid a_1\mid\cdots\mid a_j)
 =u m g(u^{-1})\mid u a_1u^{-1}\mid\cdots\mid u a_ju^{-1}.
\end{equation}
For $b(m\mid a)=m g(a)-am$, a homotopy from $T_u$ to the identity is
\begin{align*}
 H_j(m\mid a_1\mid\cdots\mid a_j)
 =\sum_{i=0}^j(-1)^{i+1}
 m g(u^{-1})
 &\mid u a_1u^{-1}\mid\cdots\mid u a_i u^{-1}\mid u\\
 &\mid a_{i+1}\mid\cdots\mid a_j.
\end{align*}
The bar differential gives $bH+Hb=T_u-\id$ by cancellation of adjacent
terms. If $u$ is a frozen monomial, this homotopy has coefficients of
nonnegative $\hbar$-adic valuation in the bases under consideration.

\subsection{A Laurent generalized Weyl algebra}

Consider
\begin{equation}\label{qt:gwa}
 B(Q,m)=F\langle z^{\pm1},x,y\rangle
 \big/\bigl(xz-Qzx,\;zy-Qyz,\;yx-1-z^m,\;
 xy-1-Q^m z^m\bigr),
\end{equation}
where $m\ge1$ and $Q=t^c$ with $c\in\Z\setminus\{0\}$.
The same argument permits a finite root of $t$, with the base field
correspondingly extended. Under the chosen formal embedding,
$Q=Q_0\exp(c\hbar)$, with $c\ne0$. Give this algebra the mutable
grading
\[
 \deg y=1,\qquad \deg x=-1,\qquad \deg z=0.
\]

\begin{lemma}\label{qt:rank-one}
The comparison in Theorem~\ref{qt:comparison} is a quasi-isomorphism
for $B(Q,m)$. Its Hochschild homology has dimensions
\[
 \dim_F\HH_0(B(Q,m)/F)=m,
 \qquad
 \dim_F\HH_1(B(Q,m)/F)
 =\dim_F\HH_2(B(Q,m)/F)=1,
\]
and vanishes in higher degrees.
\end{lemma}

\begin{proof}
Conjugation by $z$ acts on mutable weight $r$ in the bar complex by
$Q^r$. If $r\ne0$, divide the homotopy
\eqref{qt:inner-action} by $Q^r-1$. For every nonzero integer $r$,
\begin{equation}\label{qt:simple-pole}
 v_\hbar\bigl((Q^r-1)^{-1}\bigr)\ge-1.
\end{equation}
Indeed, if the constant term of $Q^r-1$ vanishes, its derivative at
$\hbar=0$ is $cr\ne0$. The resulting homotopies therefore contract
all nonzero mutable weights simultaneously, both algebraically and
continuously. It remains to compare weight zero.

Put $R=F[z^{\pm1}]$, let $\sigma(z)=Qz$, and set $a=1+z^m$.
We use the bimodule resolution and Hochschild double complex of
Solotar--Su\'arez-\'Alvarez--Vivas
\cite[Sections~4--5]{SSAV}, localized in $z$.

Let $\mathcal B$ denote the algebra in \eqref{qt:gwa} over
$R_0=E[t^{\pm1}]$, so that
\[
 B(Q,m)=\mathcal B\otimes_{R_0}F.
\]
The construction and exactness proof of
\cite[Sections~4.1--4.3]{SSAV} work over $R_0$: the filtered
Smith-algebra contraction has Laurent-polynomial coefficients in
$Q$, the central element $YX-a(z)$ is a nonzerodivisor by PBW, and
the induced vertical map on row homology is the identity. None of
these steps divides by $Q^r-1$ or by a coefficient of
$\sigma(a)-a$. Thus the same resolution remains valid at $Q=1$
and for $m=1$, although these cases are excluded from the standing
assumptions of \cite{SSAV}. Ore localization in the normal element
$z$ is exact on both sides.

Write $P_\bullet$ for the resulting free bimodule resolution and
$\operatorname{Bar}_\bullet$ for the relative bar resolution.
The PBW basis makes $\mathcal B$ free over $R_0$, so both resolutions
have free terms over
$\mathcal B\otimes_{R_0}\mathcal B^{\mathrm{op}}$.
The comparison theorem for projective resolutions, applied in the
category of mutable-weight-graded modules, therefore gives maps and
homotopies over $R_0$:
\[
 P_\bullet\mathrel{\mathop{\rightleftarrows}^{f}_{g}}
 \operatorname{Bar}_\bullet,\qquad
 gf-\id=dH+Hd,\qquad fg-\id=dJ+Jd.
\]
Each basis vector has finite image, and every scalar coefficient
$\lambda$ of $f,g,H,J$ lies in $R_0$. Under the formal embedding,
\[
 v_\hbar(\lambda)\geq0
 \qquad\text{uniformly over all basis vectors}.
\]
The maps and homotopies therefore extend to the completions, where
the resolution also computes continuous Hochschild homology.

We now use the sign convention of the displayed resolution in
\cite{SSAV}. Its vertical terms in mutable weight zero have coordinates
$R,R^3,R^3,R$. A degree-one vector is written
\[
 p_1x\mid Y+p_2\mid H+yp_3\mid X,
\]
and a degree-two vector is written
\[
 p_1x\mid Y\wedge H+p_2\mid Y\wedge X
       +yp_3\mid H\wedge X.
\]
Equations~(11)--(13) of \cite{SSAV} give
\begin{align}
 \delta_0(p)&=(\sigma(p),-a'p,\sigma(p)),\notag\\
 \delta_1(p_1,p_2,p_3)
 &=\bigl(p_1\sigma(a')+\sigma(p_2),\;
       \sigma(a)(p_3-p_1),\;
       -p_3\sigma(a')-\sigma(p_2)\bigr),
       \label{qt:vertical}\\
 \delta_2(p_1,p_2,p_3)
 &=\sigma(a)(p_1+p_3)+\sigma(a')p_2.\notag
\end{align}
Since $a'=mz^{m-1}$ is a unit, the elimination of
\cite[Proposition~5.1]{SSAV} splits~\eqref{qt:vertical}. Put
$\rho=\sigma(a)$ and $\tau=\sigma(a')$. In degree one take coordinates
\[
 u=p_1-p_3,\qquad v=p_3+\sigma(p_2/a'),\qquad w=p_2,
\]
and in degree two replace $(q_1,q_2,q_3)$ by
\[
 U=(q_1+q_3)/\tau,\qquad V=-q_3/\tau,\qquad
 W=q_2+(\rho/\tau)(q_1+q_3).
\]
The three maps become
\[
 p\longmapsto(0,0,-a'p),\qquad
 (u,v,w)\longmapsto(u,v,0),\qquad
 (U,V,W)\longmapsto\tau W.
\]
These changes and their inverses use only Laurent multiplication,
$\sigma^{\pm1}$, and inverses of $a'$ and its translates. The
resulting splittings have nonnegative $\hbar$-adic valuation.
Quotienting the top row by the vertical images
leaves a three-term complex with coordinates $R,R^2,R$ in homological
degrees $2,1,0$. To check this operation after completion, index a
term of the total complex by its exterior degree $k$ and an integer
$j\ge0$, with total degree $k+2j$ and $0\le k\le3$.
Each total degree has finitely many terms. The horizontal
differential lowers $k$ by one and preserves $j$; the vertical
differential sends $(k,j)$ to $(k+1,j-1)$. On the kernel of the
quotient map, the vertical splittings give a contraction $s$ sending
$(k,j)$ to $(k-1,j+1)$. If $d$ denotes the horizontal differential,
both $sd$ and $ds$ lower $k$ by two, so their squares vanish.
The finite perturbation formula contracts the kernel with coefficients
of nonnegative valuation, proving that the quotient is a quasi-isomorphism
before and after completion.

In degree one, use $\delta_0$ to set $p_2=0$. The horizontal
differential is
\[
 (p_1,p_3)\longmapsto D(p_1-p_3),
 \qquad D=a\sigma^{-1}-\sigma(a),
\]
where $\sigma(a)$ in the last expression acts by multiplication.
In degree two, use $p\mid Y\wedge X$ as representative. Equation~(8e)
of \cite{SSAV} gives the two coefficients $p-\sigma(p)$ and the
$H$-coefficient $-(\sigma(a)-a)'p$. Removing the $H$-coefficient
using $\delta_0$ changes each of the other two coefficients to
\[
 p-Q\frac{\sigma^2(a')}{\sigma(a')}\sigma(p)
   =(1-Q^m\sigma)(p).
\]
Thus the coordinates $(p_1-p_3,p_3)$ identify the quotient with
\begin{equation}\label{qt:two-complexes}
 \bigl[R\xrightarrow{\,1-Q^m\sigma\,}R\bigr]_{2\to1}
 \ \oplus\
 \bigl[R\xrightarrow{\,D\,}R\bigr]_{1\to0}.
\end{equation}

The first operator is diagonal on the Laurent basis:
\[
 (1-Q^m\sigma)(z^k)=(1-Q^{m+k})z^k.
\]
Its kernel and cokernel are each represented by $z^{-m}$. On all
other basis vectors its inverse has valuation at least $-1$, by
\eqref{qt:simple-pole}. Projection to $Fz^{-m}$ and this inverse give
the same retraction after completion.

For the second operator, change the source coordinate by $\sigma$.
The resulting operator is
\[
 D'=D\sigma=(1-\sigma)\circ(a\cdot).
\]
Write $e_k=z^k$ and $A_k=1-Q^k$. Then
\begin{equation}\label{qt:recurrence}
 D'e_k=A_ke_k+A_{k+m}e_{k+m}.
\end{equation}
Let $V$ be the span of $e_0,\ldots,e_{m-1}$. Define $P:R\to V$
by
\begin{equation}\label{qt:projection}
 P(e_{r+\nu m})=(-1)^\nu\frac{A_r}{A_{r+\nu m}}e_r
 \quad(1\le r<m,\ \nu\in\Z),
 \qquad
 P(e_0)=e_0,
 \qquad P(e_{\nu m})=0\quad(\nu\ne0).
\end{equation}
For $1\le r<m$ define
\begin{equation}\label{qt:homotopy-nonzero-residue}
 H(e_{r+\nu m})=
 \begin{cases}
 \displaystyle\frac{1}{A_{r+\nu m}}
   \sum_{j=0}^{\nu-1}(-1)^{\nu-1-j}e_{r+jm},&\nu>0,\\[6pt]
 0,&\nu=0,\\[3pt]
 \displaystyle\frac{1}{A_{r+\nu m}}
   \sum_{j=\nu}^{-1}(-1)^{j-\nu}e_{r+jm},&\nu<0.
 \end{cases}
\end{equation}
For the zero residue class put $H(e_0)=0$ and
\begin{equation}\label{qt:homotopy-zero-residue}
 H(e_{\nu m})=
 \begin{cases}
 \displaystyle\frac{1}{A_{\nu m}}
   \sum_{j=0}^{\nu-1}(-1)^{\nu-1-j}e_{jm},&\nu>0,\\[6pt]
 \displaystyle\frac{1}{A_{\nu m}}
   \sum_{j=\nu}^{-1}(-1)^{j-\nu}e_{jm},&\nu<0.
 \end{cases}
\end{equation}
Every sum is finite, and each denominator $A_j$ has $j\ne0$. Telescoping
\eqref{qt:recurrence} gives
\begin{equation}\label{qt:contraction}
 D'H=\id-\iota P,
 \qquad PD'=0,
 \qquad P\iota=\id_V,
\end{equation}
where $\iota:V\hookrightarrow R$. In the negative nonzero residue
class, the remaining endpoint is
$(-1)^{-1-\nu}A_rA_{r+\nu m}^{-1}e_r=-P(e_{r+\nu m})$.
For the negative zero class that endpoint vanishes because $A_0=0$.

All coefficients of $H$ and $P$ have valuation at least $-1$.
They extend to the completion, on which $D'$ is injective: if
$f=\sum_k b_ke_k$ lies in its kernel, then
\[
 b_k+b_{k-m}=0\qquad(k\ne0).
\]
In a nonzero residue class, a nonzero coefficient forces an infinite
tail of coefficients with the same valuation, contrary to
\eqref{qt:restricted-support}. In the zero residue class the positive
tail starting at $b_0$ and the negative tail ending at $b_{-m}$ must
separately vanish. Hence $f=0$. Equations~\eqref{qt:contraction} now
also imply $HD'=\id$. Thus the second complex in
\eqref{qt:two-complexes} retracts onto $V$ in degree zero, both
algebraically and continuously. The retractions commute with the
canonical inclusion of algebraic chains into completed chains, proving
the lemma.
\end{proof}

\subsection{Quantum tori and tensor products}

\begin{lemma}\label{qt:torus-products}
The comparison in Theorem~\ref{qt:comparison} holds for a
nondegenerate quantum torus, and for any finite tensor product of such
a torus and the algebras \eqref{qt:gwa}.
\end{lemma}

\begin{proof}
Let $W$ be the lattice of the quantum torus. We use its Laurent
Koszul resolution, as in \cite{Wambst}, which decomposes by total
monomial weight $v\in W$. After
normalizing the exterior generators, the summand of weight $v$ is
the Koszul complex on the scalars
\[
 t^{2\Lambda(w_i,v)}-1,
\]
where $w_i$ is a basis of $W$. Equivalently, this is the Koszul
resolution for the commuting conjugations by $Z^{w_i}$, acting on the
adjoint module of the twisted group algebra of $W$.
For $v\ne0$, nondegeneracy supplies an $i$ with
$\Lambda(w_i,v)\ne0$. Dividing by the corresponding scalar gives a
Koszul contraction with valuation at least $-1$, uniformly in $v$.
For $v=0$, the differential is zero and the complex is the finite
exterior space $\bigwedge^*W_F$. The resolution and its comparison
with bar chains are defined over the Laurent parameter ring, so the
same argument as in Lemma~\ref{qt:rank-one} proves the comparison
after completion.

For a finite tensor product, tensor the displayed retractions and use
the Hochschild Eilenberg--Zilber maps. These maps and their homotopies
are finite sums in each degree and have coefficients of nonnegative
valuation. The tensor homotopies have a uniform finite pole bound in
each degree, since the number of factors is fixed. They therefore
extend to the completed tensor products and prove the assertion.
\end{proof}

\subsection{Finite covers in the frozen directions}

We first obtain the tensor-product presentation needed in
Lemma~\ref{qt:torus-products}. Put $F_f=L_f\otimes\Q$ and
$C_f=\operatorname{span}_{\Q}(c_1,\ldots,c_n)$. Equation
\eqref{qt:compatibility} and nondegeneracy of $\Lambda$ imply
\[
 \operatorname{rad}(\Lambda|_{F_f})=C_f.
\]
Indeed, pairing with $e_1,\ldots,e_n$ injects the radical into
$\Q^n$, while it already contains the $n$-dimensional subspace
$C_f$. Choose a rational complement $W$ on which $\Lambda$ is
nondegenerate. For each $i$, add a vector of $W$ to $e_i$ so that
the resulting vector $e_i'$ is orthogonal to $W$. Write
$M_{ij}=\Lambda(e_i',e_j')$ and set
\begin{equation}\label{qt:darboux}
 \widetilde e_i=e_i'-\sum_j\frac{M_{ij}}{2d_j}c_j.
\end{equation}
Then the $\widetilde e_i$ are mutually orthogonal, remain orthogonal
to $W$, and satisfy
$\Lambda(c_i,\widetilde e_j)=d_i\delta_{ij}$.

Choose a finite-index superlattice $L_f'\supset L_f$ containing all
$\widetilde e_i-e_i$ and splitting as the direct sum of lattices on
the lines $\Q c_i$ and on $W$. A sufficiently fine common
denominator in a rational basis adapted to these subspaces gives such
a lattice. Write its basis on the $i$th coefficient line as
$c_i/m_i$, with $m_i>0$. After a finite extension of $E$, enlarge the quantum parameter by a
further finite root if necessary. Replacing $x_i$ by its product with the
frozen unit of exponent $\widetilde e_i-e_i$, and replacing $y_i$
by the inverse adjustment, yields a factor \eqref{qt:gwa} with
\[
 z_i^{m_i}=\eta_i Z^{c_i},
 \qquad Q_i=t^{-2d_i/m_i}.
\]
Scalar powers of $t$ account for the normalization of ordered
monomials. They are units at the chosen formal point. The factors
commute by \eqref{qt:darboux}, and the remaining factor is the
quantum torus on $W\cap L_f'$. Thus the enlarged algebra is a tensor
product as in Lemma~\ref{qt:torus-products}. Notice that
$-2d_i/m_i\ne0$; every rank-one parameter has nonzero logarithmic
derivative with respect to $\hbar$.

The comparison descends along finite scalar extensions $E'/E$:
both chain complexes commute with base change, since $E'[[\hbar]]$
is finite free over $E[[\hbar]]$, and $E'((\hbar))/E((\hbar))$ is
faithfully flat. We may therefore adjoin the required constant roots
and roots of unity.

Adjoining a finite root of the parameter preserves the comparison:
if $s^a=t$, choose $s_0^a=t_0$ and use
$s=s_0\exp(\hbar/a)$ in $K$. Base change of the algebraic chains
and this extension of the formal embedding give the original
comparison after tensoring with $K$.

\begin{lemma}\label{qt:frozen-descent}
Let $B$ be the algebra obtained by the above finite enlargement of the
frozen lattice. If the comparison holds for $B$, then it holds for
$A$.
\end{lemma}

\begin{proof}
Put $\Delta=L_f'/L_f$ and $G=\operatorname{Hom}(\Delta,E^\times)$.
The algebra $B$ is strongly $\Delta$-graded, its degree-zero algebra
is $A$, and each graded component contains an invertible frozen
monomial. The dual action of $G$ satisfies $B^G=A$.

In the crossed product $B\rtimes G$, the idempotent
\[
 e=\frac1{|G|}\sum_{g\in G}g
\]
is full and has corner $e(B\rtimes G)e\cong A$. To see fullness,
conjugate $e$ by frozen units representing the elements of $\Delta$.
This gives all character idempotents of $E[G]$, whose sum is one.
The resulting Morita equivalence uses finitely many units and finite
matrices, and therefore also holds after the prescribed completion.
Its maps on Hochschild complexes commute with the completion maps.

The finite group bar decomposition \cite{Lorenz} gives a quasi-isomorphism
\begin{equation}\label{qt:crossed-product}
 C_\bullet(B\rtimes G/F)
 \simeq
 \left(\bigoplus_{g\in G}C_\bullet(B,B_g)\right)^G.
\end{equation}
Here conjugation acts on the right. Group crossed-product bar chains
by their total group element and contract positive group-bar degrees
by averaging. The maps are finite sums degreewise, so the decomposition
survives completion.

Fix $g\ne1$. Some frozen unit $u$ satisfies $g(u)=\beta u$ with
$\beta\ne1$. Since $u$ is frozen, its conjugation is diagonal on
the ordered monomial basis, and the exchange relations preserve its
eigenweights by $\Lambda(L_f',c_i)=0$. On a tensor of total
eigenweight $\ell$, the operator \eqref{qt:inner-action} is therefore
multiplication by
\[
 \beta^{-1}t^\ell.
\]
The homotopy preserves these eigenweights. The function
$\beta^{-1}t^\ell-1$ is never identically zero: for
$\ell=0$ it is a nonzero constant, and for $\ell\ne0$ it is
nonconstant. At $t=t_0\exp(\hbar)$, a zero of this function is
simple, with derivative $\ell$. Its inverse consequently has
valuation at least $-1$, uniformly in the eigenweight. Dividing the
homotopy \eqref{qt:inner-action} by this function contracts the entire
$g$-summand of \eqref{qt:crossed-product}, both algebraically and
continuously.

Only the identity summand remains. Morita invariance and
\eqref{qt:crossed-product} thus identify the two sides of the
comparison for $A$ with the $G$-invariants of the corresponding sides
for $B$. The comparison for $B$ is canonical and $G$-equivariant.
Exactness of finite-group invariants in characteristic zero completes
the proof.
\end{proof}

\begin{proof}[Proof of Theorem~\ref{qt:comparison}]
The lattice enlargement above gives the tensor product of
Lemma~\ref{qt:torus-products}. Its comparison descends to $A$ by
Lemma~\ref{qt:frozen-descent}.
\end{proof}

\begingroup
\renewcommand{\bibliofont}{\footnotesize}
\raggedright\raggedbottom
\bibliographystyle{alpha}
\bibliography{references}
\endgroup
\end{document}